\documentclass{article}%
\usepackage{amsfonts}
\usepackage{amsmath}
\usepackage{amssymb}
\usepackage{graphicx}%
\newtheorem{theorem}{Theorem}

\newtheorem{example}[theorem]{Example}

\newtheorem{lemma}[theorem]{Lemma}

\newtheorem{proposition}[theorem]{Proposition}

\newenvironment{proof}[1][Proof]{\noindent\textbf{#1.} }{\ \rule{0.5em}{0.5em}}
\numberwithin{equation}{section}
\numberwithin{theorem}{section}
\begin{document}

\title{Fundamental Solutions and Mapping Properties of Semielliptic Operators}
\author{G. N. Hile\\University of Hawaii, Honolulu, HI 96822, USA\\Email : hile@hawaii.edu\\Tel :\ 808-956-6533\\Fax : 808-956-9139}
\date{}
\maketitle

\begin{abstract}
An explicit formula is given for a fundamental solution for a class of
semielliptic operators. The fundamental solution is used to investigate
properties of these operators as mappings between weighted function spaces in
$\mathbb{R}^{n}$. Necessary and sufficient conditions are given for such a
mapping to be an isomorphism. Results apply, for example, to elliptic,
parabolic, and generalized p-parabolic operators.

\end{abstract}

\noindent\textbf{Keywords:}\ Semielliptic Operator, Fundamental Solution, Isomorphism

\section{Introduction}

Let $\mathcal{L}$ denote a linear \textit{semielliptic }partial differential
operator, acting on suitable real or complex $m\times1$ vector functions
$u=u\left(  x\right)  $, $x\in\mathbb{R}^{n}$, according to%
\begin{equation}
\mathcal{L}u=\sum_{\alpha\cdot\gamma=\ell}A_{\alpha}\partial^{\alpha}%
u\quad.\label{op}%
\end{equation}
The coefficients $\left\{  A_{\alpha}\right\}  $ are constant $m\times m$
matrices with real or complex entries, indexed by multi-indices $\alpha
=\left(  \alpha_{1},\cdots,\alpha_{n}\right)  $ in $\mathbb{R}^{n}$. The
positive integer $\ell$\ is the \textit{order} of $\mathcal{L}$,
\[
\ell=\max\left\{  \left\vert \alpha\right\vert :A_{\alpha}\neq0\right\}
\quad,
\]
and $\gamma$ is a fixed vector of rational numbers,
\begin{equation}
\gamma=\left(  \gamma_{1},\gamma_{2},\cdots,\gamma_{n}\right)  =\left(
\frac{\ell}{\ell_{1}},\frac{\ell}{\ell_{2}},\ldots,\frac{\ell}{\ell_{n}%
}\right)  \quad,\label{gamma}%
\end{equation}
with each $\ell_{k}$ a positive integer. The \textit{semiellipticity
condition} on $\mathcal{L}$ requires that its \textit{symbol}, the matrix of
polynomial functions%
\begin{equation}
L\left(  x\right)  =\sum_{\alpha\cdot\gamma=\ell}A_{\alpha}\left(  ix\right)
^{\alpha}\quad,\label{ell}%
\end{equation}
be invertible for all nonzero $x$ in $\mathbb{R}^{n}$. (Alternative adjectives
to \textit{semielliptic}, all used by various authors, are
\textit{quasielliptic}, \textit{semi-elliptic}, and \textit{quasi-elliptic}).
As explained for example in \cite{HMZ}, \S 2, a consequence of the
semiellipticity condition is that $\max_{k}\ell_{k}=\ell$, so that $\gamma
_{k}\geq1$ for each $k$ and $\gamma_{k}=1$ for at least one $k$. A further
consequence is that, for each $k$, $1\leq k\leq n$, the term%
\[
A_{\ell_{k}e_{k}}\partial^{\ell_{k}}/\partial x_{k}^{\ell_{k}}\qquad,
\]
corresponding to $\alpha=\ell_{k}e_{k}$ where $e_{k}$ is the $k$th unit
coordinate vector in $\mathbb{R}^{n}$, appears in $\mathcal{L}$ with
$A_{\ell_{k}e_{k}}$ an invertible matrix. This term is the only unmixed
derivative with respect to $x_{k}$ appearing in $\mathcal{L}$, and $\ell_{k} $
is the highest order of differentiation with respect to $x_{k}$ appearing in
$\mathcal{L}$.

We consider the subclass of such operators satisfying the additional condition%
\begin{equation}
\left\Vert \gamma\right\Vert :=\sum_{k=1}^{n}\gamma_{k}>\ell\quad
.\label{gamcond}%
\end{equation}
We show that for such operators $\mathcal{L}$ a fundamental solution is given
explicitly by the iterated integral%
\begin{equation}
F\left(  x\right)  =\left(  2\pi\right)  ^{-n}\int_{0}^{\infty}\int
_{\mathbb{R}^{n}}e^{ix\cdot z}\sigma\left(  z\right)  e^{-t\sigma\left(
z\right)  }L\left(  z\right)  ^{-1}\ dz\ dt\qquad\left(  x\neq0\right)
\quad,\label{eff}%
\end{equation}
where $\sigma$ is the function%
\begin{equation}
\sigma\left(  x\right)  =\sum_{k=1}^{n}x_{k}{}^{2\ell_{k}}\quad.\label{sigma}%
\end{equation}
In particular, $F$ has the properties

\begin{itemize}
\item $F\in C^{\infty}\left(  \mathbb{R}^{n}\backslash\left\{  0\right\}
\right)  $,

\item $\mathcal{L}F\left(  x\right)  =0$ for $x\neq0$,

\item for all complex $m\times1$ vector functions $\varphi$ in $C_{0}^{\infty
}\left(  \mathbb{R}^{n}\right)  $, $F\ast\varphi\in C^{\infty}\left(
\mathbb{R}^{n}\right)  $ and $\mathcal{L}\left(  F\ast\varphi\right)
=\varphi$.
\end{itemize}

\noindent It must be pointed out that the order of integration is important in
(\ref{eff}), as Fubini's Theorem does not apply, and interchanging orders of
integration will likely destroy convergence of the integrals.

Of course other methods are known for constructing fundamental solutions for
partial differential equations with constant coefficients. Unfortunately these
methods do not always produce representations highly useful for investigation
of solutions of the corresponding equations.

We use our fundamental solution to investigate properties of the operator
$\mathcal{L}$ as a mapping between weighted function spaces. We introduce the
vector of integers%
\begin{equation}
\underline{\ell}=\left(  \ell_{1},\ell_{2},\ldots,\ell_{n}\right)
\quad,\label{ellvec}%
\end{equation}
and a corresponding weight function%
\begin{equation}
\rho\left(  x\right)  =\sigma\left(  x\right)  ^{1/\left(  2\ell\right)
}=\left(  \sum_{k=1}^{n}x_{k}{}^{2\ell_{k}}\right)  ^{1/\left(  2\ell\right)
}\quad.\label{rho}%
\end{equation}
We define function spaces $W_{s}^{r,p}\left(  \mathbb{R}^{n},\mathbb{C}%
^{m},\underline{\ell}\right)  $, consisting of complex $m\times1$ vector
functions $u$ on $\mathbb{R}^{n}$ with finite norm%
\begin{equation}
\left\Vert u\right\Vert _{r,p,s;\underline{\ell}}=\sum_{\alpha\cdot\gamma\leq
r}\left\Vert \left(  1+\rho\right)  ^{s+\alpha\cdot\gamma}\partial^{\alpha
}u\right\Vert _{p,\mathbb{R}^{n}}\quad.\label{norm}%
\end{equation}
We demonstrate (Theorem \ref{main}) that, if $1<p<\infty$ and $\left\Vert
\gamma\right\Vert >\ell$, then the mapping%
\begin{equation}
\mathcal{L}:W_{s}^{\ell,p}\left(  \mathbb{R}^{n},\mathbb{C}^{m},\underline
{\ell}\right)  \longrightarrow W_{s+\ell}^{0,p}\left(  \mathbb{R}%
^{n},\mathbb{C}^{m},\underline{\ell}\right) \label{mapit}%
\end{equation}
is an isomorphism if and only if
\begin{equation}
-\left\Vert \gamma\right\Vert /p<s<\left\Vert \gamma\right\Vert -\ell
-\left\Vert \gamma\right\Vert /p\quad.\label{conit}%
\end{equation}

Investigations of mappings between weighted Sobolev spaces, analogous to
(\ref{mapit}), began with Cantor \cite{CA}, who considered the scalar operator
$\mathcal{L}=\Delta$, the Laplace operator in $\mathbb{R}^{n}$. Then $m=1$,
$\ell=2$, $\gamma=\left(  1,\ldots,1\right)  $, and the weight function $\rho$
is equivalent to the Euclidean norm. Cantor showed that in this special case
(\ref{mapit}) is an isomorphism provided that $n>2$, $n/\left(  n-2\right)
<p<\infty$, and $-n/p<s<n-2-n/p$. Cantor was building on the work of Walker
and Nirenberg \cite{WA1,WA2,NW}, who showed that certain elliptic differential
operators have finite dimensional null spaces as mappings between (unweighted)
Sobolev spaces $W^{\ell,p}\left(  \mathbb{R}^{n}\right)  $.

McOwen \cite{MC1} removed Cantor's restriction $p>n/\left(  n-2\right)  $ for
the Laplace operator, and further listed various conditions, similar to
(\ref{conit}), that guarantee a power $\Delta^{k}$ of the Laplacian is a
Fredholm map having certain properties. In particular, for $\Delta^{k}$ we
have $\underline{\ell}=\left(  2k,\ldots,2k\right)  $, and if $n>2k$ and
$-n/p<s<n-2k-n/p$, the map%
\[
\Delta^{k}:W_{s}^{2k,p}\left(  \mathbb{R}^{n},\mathbb{C},\underline{\ell
}\right)  \longrightarrow W_{s+2k}^{0,p}\left(  \mathbb{R}^{n},\mathbb{C}%
,\underline{\ell}\right)
\]
is an isomorphism. Lockhart \cite{LOC} extended this result to scalar elliptic
operators of order $\ell$, with constant coefficients and only highest order
terms; for these operators (\ref{mapit}) is an isomorphism provided that
$-n/p<s<n-\ell-n/p$. Lockhart and McOwen \cite{LOC,LM,MC2} consider also
elliptic operators with variable coefficients continuous at infinity; for such
operators conditions are given under which the mapping (\ref{mapit}) is Fredholm.

In \cite{HM}, the author and C. Mawata considered the mapping (\ref{mapit})
for the case of the heat operator in $\mathbb{R}^{n}$, $\mathcal{L}%
=\partial_{t}-\Delta u$. They gave conditions under which the mapping is
Fredholm, and in particular showed that $\mathcal{L}$ is an isomorphism
provided that $-\left(  n+2\right)  /p<s<n-\left(  n+2\right)  /p$.

In the last section of the paper are examples demonstrating how the main
isomorphism theorem extends known results on elliptic operators, and produces
new results for parabolic and generalized r-parabolic operators.

Of special relevance to this work is that of G. V. Demidenko \cite{DEM2,DEM3}
who, working in weighted Sobolev spaces somewhat different from ours, studied
mapping properties of the semielliptic opertor (\ref{op}). When translated to
the notation of this paper, his result on isomorphic properties asserts that
the mapping%
\[
\mathcal{L}:W_{-\ell}^{\ell,p}\left(  \mathbb{R}^{n},\mathbb{C}^{m}%
,\underline{\ell}\right)  \longrightarrow L^{p}\left(  \mathbb{R}%
^{n},\mathbb{C}^{m}\right)
\]
is an isomorphism provided that $1<p<\left\Vert \gamma\right\Vert /\ell$. This
statement is a special case of our isomorphism result, obtained by taking
$s=-\ell$ in (\ref{mapit}) and (\ref{conit}). Demidenko did not use
fundamental solutions in his investigations, but rather used integral
representations to present what he called \textquotedblleft approximate
solutions\textquotedblright\ of equations $\mathcal{L}u=f$, converging in
$L^{p}$ to actual solutions. A modification of Demidenko's representations led
to our discovery of formula (\ref{eff}) for a fundamental solution.

Demidenko \cite{DEM4,DEM5,DEM6} has recently extended his isomorphism results
to operators of the form%
\[
\mathcal{L}=\left[
\begin{tabular}
[c]{cc}%
$\mathcal{L}_{1}$ & $\mathcal{L}_{2}$\\
$\mathcal{L}_{3}$ & $0$%
\end{tabular}
\ \ \right]  \qquad,
\]
where $\mathcal{L}_{1}$ is a square matrix semielliptic operator, and
$\mathcal{L}_{2}$ and $\mathcal{L}_{3}$ are rectangular matrix differential
operators having certain properties related to semiellipticity.

\section{Notation and Preliminaries}

In formula (\ref{op}) we use the conventional notation%
\[
\alpha\cdot\gamma=\alpha_{1}\gamma_{1}+\alpha_{2}\gamma_{2}+\cdots+\alpha
_{n}\gamma_{n}\qquad,\qquad\partial^{\alpha}=\frac{\partial^{\left\vert
\alpha\right\vert }}{\partial x_{1}^{\alpha_{1}}\partial x_{2}^{\alpha_{2}%
}\cdots\partial x_{n}^{\alpha_{n}}}\quad,
\]
with $\left\vert \alpha\right\vert =\alpha_{1}+\alpha_{2}+\cdots+\alpha_{n}$
representing the \textit{length} of the multi-index $\alpha$. Some authors
write (\ref{op}) in the equivalent formulation%
\[
\mathcal{L}u=\sum_{\alpha/\underline{\ell}=1}A_{\alpha}\partial^{\alpha}%
u\quad,
\]
where $\alpha/\underline{\ell}$ is defined as the sum $\alpha_{1}/\ell
_{1}+\cdots+\alpha_{n}/\ell_{n}$. However we prefer (\ref{op}), as the vector
$\gamma$ proves useful in some of our representations.

For a vector $x\in\mathbb{R}^{n}$ and for a complex matrix $M=\left(
m_{ij}\right)  $, we employ the usual norms%
\[
\left\vert x\right\vert =\left(  x_{1}^{2}+x_{2}^{2}+\cdots+x_{n}^{2}\right)
^{1/2}\qquad,\qquad\left\vert M\right\vert =\left(  \sum_{i,j}\left\vert
m_{ij}\right\vert ^{2}\right)  ^{1/2}\quad.
\]
We also at times use an alternate norm for vectors, as specified by%
\[
\left\Vert x\right\Vert :=\left\vert x_{1}\right\vert +\left\vert
x_{2}\right\vert +\cdots+\left\vert x_{n}\right\vert \quad,
\]
and as demonstrated already in (\ref{gamcond}). (The length $\left\vert
\alpha\right\vert $ of a multi-index $\alpha$ is not the Euclidean length, but
rather the same as $\left\Vert \alpha\right\Vert $; however we conform to
custom and use $\left\vert \alpha\right\vert $, with the expectation that the
correct interpretation will be clear from the context.)

The function $\rho$ of (\ref{rho}) serves as an \textit{anisotropic length} of
vectors $x$ in $\mathbb{R}^{n}$. In \cite{HMZ}, \S 2 and \S 5, one finds
verifications of the inequalities%
\begin{equation}
\rho(x+y)\leq\rho(x)+\rho(y)\qquad,\qquad\left\vert x^{\alpha}\right\vert
\leq\rho(x)^{\alpha\cdot\gamma}\quad.\label{tri}%
\end{equation}
For positive real numbers $t$ and for $x\in\mathbb{R}^{n}$ we denote%
\[
t^{\gamma}x:=\left(  t^{\gamma_{1}}x_{1},t^{\gamma_{2}}x_{2},\cdots
,t^{\gamma_{n}}x_{n}\right)  \quad.
\]
Straightforward calculations confirm that%
\begin{equation}
\left(  t^{\gamma}x\right)  ^{\alpha}=t^{\alpha\cdot\gamma}x^{\alpha}%
\qquad,\qquad\rho\left(  t^{\gamma}x\right)  =t\rho(x)\qquad,\qquad L\left(
t^{\gamma}x\right)  =t^{\ell}L\left(  x\right)  \quad.\label{hom}%
\end{equation}
If we fix $x$ in the last two of these inequalities and choose $t=1/\rho
\left(  x\right)  $, so that $\rho\left(  t^{\gamma}x\right)  =1$, we deduce
that
\begin{equation}
c_{1}\left(  \mathcal{L}\right)  \rho(x)^{\ell}\leq\left\vert L\left(
x\right)  \right\vert \leq c_{2}\left(  \mathcal{L}\right)  \rho(x)^{\ell
}\quad,\label{rhoineq1}%
\end{equation}%
\begin{equation}
c_{3}\left(  \mathcal{L}\right)  \rho(x)^{-\ell}\leq\left\vert L\left(
x\right)  ^{-1}\right\vert \leq c_{4}\left(  \mathcal{L}\right)
\rho(x)^{-\ell}\qquad\left(  x\neq0\right)  \quad,\label{rhoineq2}%
\end{equation}
where $\left\vert \cdots\right\vert $ here is the matrix norm, and
\begin{align*}
c_{1}\left(  \mathcal{L}\right)   &  =\min_{\rho(y)=1}\ \left\vert
L(y)\right\vert \qquad,\qquad c_{2}\left(  \mathcal{L}\right)  =\max
_{\rho(y)=1}\ \left\vert L(y)\right\vert \quad,\\
c_{3}\left(  \mathcal{L}\right)   &  =\min_{\rho(y)=1}\ \left\vert
L(y)^{-1}\right\vert \qquad,\qquad c_{4}\left(  \mathcal{L}\right)
=\max_{\rho(y)=1}\ \left\vert L(y)^{-1}\right\vert \quad.
\end{align*}
It is further demonstrated in \cite{HMZ}, \S 2, that, for $x$ in
$\mathbb{R}^{n}$ and multi-indices $\alpha$, there are nonnegative constants
$c_{5}\left(  \mathcal{L}\right)  $ and $c_{6}\left(  \mathcal{L}%
,\alpha\right)  $ such that%
\begin{equation}
\left\vert \partial^{\alpha}L(x)\right\vert \leq\left\{
\begin{tabular}
[c]{cc}%
$c_{5}\left(  \mathcal{L}\right)  \rho\left(  x\right)  ^{\ell-\alpha
\cdot\gamma}$ & if $\alpha\cdot\gamma\leq\ell,$\\
$0$ & otherwise,
\end{tabular}
\right. \label{rhoineq3}%
\end{equation}%
\begin{equation}
\left\vert \partial^{\alpha}L\left(  x\right)  ^{-1}\right\vert \leq
c_{6}\left(  \mathcal{L},\alpha\right)  \rho\left(  x\right)  ^{-\ell
-\alpha\cdot\gamma}\qquad\left(  x\neq0\right)  \quad.\label{rhoineq4}%
\end{equation}

\section{A Fundamental Solution}

We demonstrate that formula (\ref{eff}) does indeed prescribe a fundamental
solution for the differential operator $\mathcal{L}$. First we investigage in
some detail the inner integral of (\ref{eff}). For multi-indices $\beta$ and
points $x$ in $\mathbb{R}^{n}$, and for $t>0$, we define $m\times m$ matrix
valued functions%

\begin{align}
J\left(  x,t\right)   &  =\int_{\mathbb{R}^{n}}e^{ix\cdot z}\sigma
(z)e^{-t\sigma(z)}L\left(  z\right)  ^{-1}\ dz\quad,\label{ja}\\
J_{\beta}\left(  x,t\right)   &  =\int_{\mathbb{R}^{n}}e^{ix\cdot z}\left(
iz\right)  ^{\beta}\sigma(z)e^{-t\sigma(z)}L\left(  z\right)  ^{-1}%
\ dz\quad.\label{jay}%
\end{align}
Observe that $J_{\beta}=J$ when $\beta=\left(  0,\ldots,0\right)  $. Also,
$J_{\beta}$ is the formal derivative $\partial_{x}^{\beta}J$, arising by
differentiation of $J$ under the integral; we will discuss validity of this
action, as well as convergence of these integrals and other properties.

We require a lemma concerning convergence of more elementary integrals.

\begin{lemma}
\label{rhoint}Let $s$ be any real constant, and assume $t>0$.

(a)\ The integral%
\[
\int_{\rho\left(  x\right)  \geq1}\rho\left(  x\right)  ^{s}\ dx
\]
is finite if and only if $s<-\left\Vert \gamma\right\Vert $, in which case%
\begin{equation}
\int_{\rho\left(  x\right)  \geq t}\rho\left(  x\right)  ^{s}%
\ dx=t^{s+\left\Vert \gamma\right\Vert }\int_{\rho\left(  x\right)  \geq1}%
\rho\left(  x\right)  ^{s}\ dx\quad.\label{ja2}%
\end{equation}

(b)\ The integral%
\[
\int_{\rho\left(  x\right)  \leq1}\rho\left(  x\right)  ^{s}\ dx
\]
is finite if and only if $s>-\left\Vert \gamma\right\Vert $, in which case%
\begin{equation}
\int_{\rho\left(  x\right)  \leq t}\rho\left(  x\right)  ^{s}%
\ dx=t^{s+\left\Vert \gamma\right\Vert }\int_{\rho\left(  x\right)  \leq1}%
\rho\left(  x\right)  ^{s}\ dx\quad.\label{ja1}%
\end{equation}

\end{lemma}

\begin{proof}
Given $0\leq r<R\leq\infty$ and $s\in\mathbb{R}$, let $I\left(  r,R,s\right)
$ be the integral%
\begin{equation}
I\left(  r,R,s\right)  =\int_{r\leq\rho\left(  x\right)  \leq R}\rho\left(
x\right)  ^{s}\ dx\quad.\label{ei}%
\end{equation}
In this integral we make the change of integration parameter $z=t^{\gamma}x$,
where $t$ is a positive constant; then $dz=t^{\left\Vert \gamma\right\Vert
}\ dx$, $\rho\left(  z\right)  =t\rho\left(  x\right)  $, to derive%
\begin{equation}
I\left(  tr,tR,s\right)  =t^{s+\left\Vert \gamma\right\Vert }I\left(
r,R,s\right)  \quad.\label{ja3}%
\end{equation}
Setting $r=1$, $R=2$, and $t=2^{m}$ gives%
\[
I\left(  2^{m},2^{m+1},s\right)  =2^{m\left(  s+\left\Vert \gamma\right\Vert
\right)  }I\left(  1,2,s\right)  \quad.
\]
From the geometric summation%
\[
I\left(  1,\infty,s\right)  =\sum_{m=0}^{\infty}I\left(  2^{m},2^{m+1}%
,s\right)  =I\left(  1,2,s\right)  \sum_{m=0}^{\infty}2^{m\left(  s+\left\Vert
\gamma\right\Vert \right)  }%
\]
it follows that $I\left(  1,\infty,s\right)  <\infty$ if and only if
$s+\left\Vert \gamma\right\Vert <0$. Then we obtain (\ref{ja2}) by setting
$r=1$, $R=\infty$ in (\ref{ja3}). Likewise, from%
\[
I\left(  0,1,s\right)  =\sum_{m=-1}^{-\infty}I\left(  2^{m},2^{m+1},s\right)
=I\left(  1,2,s\right)  \sum_{m=-1}^{-\infty}2^{m\left(  s+\left\Vert
\gamma\right\Vert \right)  }%
\]
it follows that $I\left(  0,1,s\right)  <\infty$ if and only if $s+\left\Vert
\gamma\right\Vert >0$, in which case (\ref{ja1}) is obtained by setting $r=0
$, $R=1$ in (\ref{ja3}).
\end{proof}

\qquad

The next lemma is proved in \cite{HMZ}, \S 5.

\begin{lemma}
\label{hmzlemma}Given real numbers R and S with $0\leq R<S$, there exists a
real valued function $\psi$ in $C_{0}^{\infty}\left(  \mathbb{R}^{n}\right)
$, with support in the region where $\rho\left(  x\right)  <S$, such that
$0\leq\psi\leq1$, $\psi(x)=1$ if $\rho(x)\leq R$, and for any multi-index
$\alpha$ and $x$ in $\mathbb{R}^{n}$,%
\[
\left\vert \partial^{\alpha}\psi(x)\right\vert \leq C(\ell,\alpha)\left(
S-R\right)  ^{-\alpha\cdot\gamma}\quad.
\]

\end{lemma}

\qquad

The following lemma, somewhat technical in nature, gathers pertinent
information regarding the integrals $J_{\beta}$.

\begin{lemma}
\label{jaylemma}Each integral $J_{\beta}\left(  x,t\right)  $ converges
absolutely for $x\in\mathbb{R}^{n}$ and $t>0$, with%
\begin{equation}
\left\vert J_{\beta}\left(  x,t\right)  \right\vert \leq C\left(
\mathcal{L},\beta\right)  \frac{t^{-1/2}}{\left[  t^{1/\left(  2\ell\right)
}+\rho\left(  x\right)  \right]  ^{\beta\cdot\gamma+\left\Vert \gamma
\right\Vert }}\quad.\label{jaybound}%
\end{equation}
Moreover, $J\in C^{\infty}\left[  \mathbb{R}^{n}\times\left(  0,\infty\right)
\right]  $, with differentiation of $J$ under the integral of all orders
allowed. In particular, for $x\in\mathbb{R}^{n}$ and $t>0$,%
\begin{align}
\frac{\partial^{\beta}}{\partial x^{\beta}}J\left(  x,t\right)   &  =J_{\beta
}\left(  x,t\right)  \quad,\label{jayy}\\
\frac{\partial^{k}}{\partial t^{k}}\frac{\partial^{\beta}}{\partial x^{\beta}%
}J\left(  x,t\right)   &  =\left(  -1\right)  ^{k}\int_{\mathbb{R}^{n}%
}e^{ix\cdot z}\left(  iz\right)  ^{\beta}\sigma\left(  z\right)
^{k+1}e^{-t\sigma\left(  z\right)  }L\left(  z\right)  ^{-1}\ dz\quad
,\label{jayyy}%
\end{align}
with (\ref{jayyy}) likewise converging absolutely. If also $s>0$, then%
\begin{equation}
J_{\beta}\left(  x,t\right)  =s^{\beta\cdot\gamma+\left\Vert \gamma\right\Vert
+\ell}J_{\beta}\left(  s^{\gamma}x,s^{2\ell}t\right)  \quad.\label{jayz}%
\end{equation}

\end{lemma}

\begin{proof}
Using (\ref{rho}), (\ref{tri}), and (\ref{rhoineq2}), we derive for the
integrand of (\ref{jayyy}) the bound%
\begin{equation}
\left\vert e^{ix\cdot z}\left(  iz\right)  ^{\beta}\sigma(z)^{k+1}%
e^{-t\sigma(z)}L\left(  z\right)  ^{-1}\right\vert \leq c_{4}\left(
\mathcal{L}\right)  \rho\left(  z\right)  ^{\beta\cdot\gamma+2k\ell+\ell
}e^{-t\sigma\left(  z\right)  }\quad.\label{jayb}%
\end{equation}
This exponential decay confirms the absolute convergence of (\ref{jayyy}), and
(when $k=0$) of $J_{\beta}\left(  x,t\right)  $.

Next we write%
\[
J_{\beta}\left(  s^{\gamma}x,s^{2\ell}t\right)  =\int_{\mathbb{R}^{n}%
}e^{is^{\gamma}x\cdot z}\left(  iz\right)  ^{\beta}\sigma(z)e^{-s^{2\ell
}t\sigma(z)}L\left(  z\right)  ^{-1}\ dz\quad,
\]
and in the integral make the change of variable $y=s^{\gamma}z$, with%
\[
y^{\beta}=s^{\beta\cdot\gamma}z^{\beta}\quad,\quad\sigma\left(  y\right)
=s^{2\ell}\sigma\left(  z\right)  \quad,\quad L\left(  y\right)  =s^{\ell
}L\left(  z\right)  \quad,\quad dy=s^{\left\Vert \gamma\right\Vert }dz\quad,
\]
to obtain (\ref{jayz}).

We consider differentiating the integral on the right of (\ref{jayyy}), which
we will refer to as $Q\left(  x,t\right)  $, with respect to $x_{k}$.
Recalling that $\left\vert e^{ir}-1\right\vert \leq\left\vert r\right\vert $
for $r\in\mathbb{R}$, we obtain for $h\neq0$ the estimate%
\[
\left\vert \frac{1}{h}\left[  Q\left(  x+he_{k},t\right)  -Q\left(
x,t\right)  \right]  \right\vert \leq\int_{\mathbb{R}^{n}}\left\vert
z_{k}\right\vert \left\vert e^{ixz}\left(  iz\right)  ^{\beta}\sigma
(z)^{k+1}e^{-t\sigma(z)}L\left(  z\right)  ^{-1}\right\vert \ dz,
\]
and note that the latter integral converges in view of (\ref{jayb}). By the
dominated convergence theorem, differentiation of $Q\left(  x,t\right)  $ with
respect to $x_{k}$ under the integral is valid. In particular,%
\[
\frac{\partial}{\partial x_{k}}J_{\beta}\left(  x,t\right)  =J_{\beta+e_{k}%
}\left(  x,t\right)  \quad,
\]
and an induction argument confirms (\ref{jayy}). In a similar way,
differentiation of $Q\left(  x,t\right)  $ under the integral with respect to
$t$ can be justified with use of the inequality $\left\vert e^{z}-1\right\vert
\leq\left\vert z\right\vert e^{\left\vert z\right\vert }$; then (\ref{jayyy})
follows by induction.

It remains only to establish the bound (\ref{jaybound}). First we bound the
integral%
\begin{equation}
J_{\beta}\left(  x,1\right)  =\int_{\mathbb{R}^{n}}e^{ix\cdot z}\left(
iz\right)  ^{\beta}\sigma(z)e^{-\sigma(z)}L\left(  z\right)  ^{-1}%
\ dz\quad.\label{jay1}%
\end{equation}

By Lemma \ref{hmzlemma}, there exists a real valued function $\psi$ in
$C_{0}^{\infty}\left(  \mathbb{R}^{n}\right)  $, with $0\leq\psi\leq1$,
$\psi\left(  z\right)  =0$ if $\rho\left(  z\right)  \geq2$, $\psi\left(
z\right)  =1$ if $\rho\left(  z\right)  \leq1$, and $\left\vert \partial
^{\alpha}\psi\left(  z\right)  \right\vert \leq C\left(  \ell,\alpha\right)  $
for any multi-index $\alpha$ in $\mathbb{R}^{n}$. Given $\varepsilon>0$, we
set
\[
\varphi_{\varepsilon}\left(  z\right)  =1-\psi\left[  \left(  \varepsilon
^{-1}\right)  ^{\gamma}z\right]  =1-\psi\left(  \frac{z_{1}}{\varepsilon
^{\gamma_{1}}},\frac{z_{2}}{\varepsilon^{\gamma_{2}}},\cdots,\frac{z_{n}%
}{\varepsilon^{\gamma_{n}}}\right)  \quad,
\]
so that $\varphi_{\varepsilon}\in C^{\infty}\left(  \mathbb{R}^{n}\right)  $,
$0\leq\varphi_{\varepsilon}\leq1$, and%
\begin{equation}
\varphi_{\varepsilon}\left(  z\right)  =\left\{
\begin{array}
[c]{cc}%
0\quad, & \text{if }\rho\left(  z\right)  \leq\varepsilon\\
1\quad, & \text{if }\rho\left(  z\right)  \geq2\varepsilon
\end{array}
\right.  \quad,\quad\left\vert \partial^{\alpha}\varphi_{\varepsilon}\left(
z\right)  \right\vert \leq C\left(  \ell,\alpha\right)  \varepsilon
^{-\alpha\cdot\gamma}\quad.\label{phi}%
\end{equation}
But if $\alpha\neq0$, $\partial^{\alpha}\varphi_{\varepsilon}(z)$ vanishes
except where $\varepsilon\leq\rho(z)\leq2\varepsilon$; thus the second
inequality of (\ref{phi}) implies also%
\begin{equation}
\left\vert \partial^{\alpha}\varphi_{\varepsilon}\left(  z\right)  \right\vert
\leq C\left(  \ell,\alpha\right)  \rho\left(  z\right)  ^{-\alpha\cdot\gamma
}\quad.\label{phib}%
\end{equation}
(For $\alpha=0$ inequality (\ref{phib}) is trivial.)

As (\ref{jay1}) converges absolutely, for any multi-index $\alpha$ we may
write%
\begin{align*}
\left(  ix\right)  ^{\alpha}J_{\beta}\left(  x,1\right)   &  =\left(
ix\right)  ^{\alpha}\ \lim_{\varepsilon\rightarrow0}\ \int_{\mathbb{R}^{n}%
}\varphi_{\varepsilon}\left(  z\right)  e^{ix\cdot z}\left(  iz\right)
^{\beta}\sigma(z)e^{-\sigma(z)}L\left(  z\right)  ^{-1}\ dz\\
&  =\lim_{\varepsilon\rightarrow0}\ \int_{\rho\left(  z\right)  \geq
\varepsilon}\left(  \partial_{z}^{\alpha}e^{ix\cdot z}\right)  \varphi
_{\varepsilon}\left(  z\right)  \left(  iz\right)  ^{\beta}\sigma
(z)e^{-\sigma(z)}L\left(  z\right)  ^{-1}\ dz\quad.
\end{align*}
We integrate by parts, taking into account the exponential decay of the
integrand at infinity (which we discuss in more detail later), as well as the
vanishing of $\varphi_{\varepsilon}\left(  z\right)  $ and all its derivatives
on the surface $\rho\left(  z\right)  =\varepsilon$, to obtain%
\begin{align}
&  \left(  ix\right)  ^{\alpha}J_{\beta}\left(  x,1\right) \nonumber\\
&  =\left(  -1\right)  ^{\left\vert \alpha\right\vert }\lim_{\varepsilon
\rightarrow0}\int_{\rho\left(  z\right)  \geq\varepsilon}e^{ix\cdot z}%
\partial^{\alpha}\left[  \varphi_{\varepsilon}\left(  z\right)  \left(
iz\right)  ^{\beta}\sigma(z)e^{-\sigma(z)}L\left(  z\right)  ^{-1}\right]
\ dz\ .\label{jbeta}%
\end{align}

Now, the derivative
\[
\partial^{\alpha}\left[  \varphi_{\varepsilon}\left(  z\right)  \left(
iz\right)  ^{\beta}\sigma(z)e^{-\sigma(z)}L\left(  z\right)  ^{-1}\right]
\]
is a finite linear combination of products of the form%
\[
\partial^{\eta}\varphi_{\varepsilon}\left(  z\right)  \ \partial^{\mu}%
z^{\beta}\ \partial^{\nu}\left[  \sigma\left(  z\right)  e^{-\sigma\left(
z\right)  }\right]  \ \partial^{\tau}L\left(  z\right)  ^{-1}\quad,
\]
where $\eta$, $\mu$, $\nu$, and $\tau$ are multi-indices in $\mathbb{R}^{n}$
with $\eta+\mu+\nu+\tau=\alpha$. From (\ref{phib}) and (\ref{rhoineq4}),%
\begin{equation}
\left\vert \partial^{\eta}\varphi_{\varepsilon}\left(  z\right)  \right\vert
\leq C\left(  \ell,\eta\right)  \rho\left(  z\right)  ^{-\eta\cdot\gamma}%
\quad,\quad\left\vert \partial^{\tau}L\left(  z\right)  ^{-1}\right\vert \leq
C\left(  \mathcal{L},\tau\right)  \rho\left(  z\right)  ^{-\ell-\tau
\cdot\gamma}\ .\label{bound1}%
\end{equation}
In view of the formula%
\[
\partial^{\mu}z^{\beta}=\left\{
\begin{array}
[c]{cc}%
\frac{\beta!}{\left(  \beta-\mu\right)  !}z^{\beta-\mu} & ,\quad\text{if }%
\mu\leq\beta\quad\text{,}\\
0 & ,\quad\text{otherwise\quad,}%
\end{array}
\right.
\]
we may use the second inequality of (\ref{tri}) to obtain
\begin{equation}
\left\vert \partial^{\mu}z^{\beta}\right\vert \leq C\left(  \beta\right)
\rho\left(  z\right)  ^{\left(  \beta-\mu\right)  \cdot\gamma}\quad
.\label{bound2}%
\end{equation}

From (\ref{rho}) and (\ref{hom}), for $t>0$ and any multi-index $\omega$ we
infer that%
\[
\sigma\left(  z\right)  =t^{-2\ell}\sigma\left(  t^{\gamma}z\right)
\qquad,\qquad\partial^{\omega}\sigma\left(  z\right)  =t^{\omega\cdot
\gamma-2\ell}\partial^{\omega}\sigma\left(  t^{\gamma}z\right)  \quad.
\]
If $z\neq0$ we may choose $t=\sigma\left(  z\right)  ^{-1/(2\ell)}=\rho\left(
z\right)  ^{-1}$, so that $\rho\left(  t^{\gamma}z\right)  =1$; then we obtain%
\begin{equation}
\left\vert \partial^{\omega}\sigma\left(  z\right)  \right\vert \leq C\left(
\omega,\underline{\ell}\right)  \rho\left(  z\right)  ^{2\ell-\omega
\cdot\gamma}\quad,\label{bound3}%
\end{equation}
where $C\left(  \omega,\underline{\ell}\right)  =\sup_{\rho\left(  z\right)
=1}\left\vert \partial^{\omega}\sigma\left(  z\right)  \right\vert $.

Any derivative $\partial^{\nu}\left[  \sigma\left(  z\right)  e^{-\sigma
\left(  z\right)  }\right]  $ is a finite linear combination of terms of the
form%
\[
e^{-\sigma\left(  z\right)  }\left[  \prod_{k=1}^{K}\partial^{\nu^{k}}%
\sigma\left(  z\right)  \right]  \quad,
\]
where $1\leq K\leq1+\left\vert \nu\right\vert \leq1+\left\vert \alpha
\right\vert $ and $\left\{  \nu^{k}\right\}  $ are multi-indices with $\nu
^{1}+\nu^{2}+\cdots+\nu^{K}=\nu$. As everything in this linear combination
depends upon $\nu$ and $\underline{\ell}$, we obtain with use of
(\ref{bound3}) the estimate%
\begin{align*}
\left\vert e^{-\sigma\left(  z\right)  }\left[  \prod_{k=1}^{K}\partial
^{\nu^{k}}\sigma\left(  z\right)  \right]  \right\vert  &  \leq e^{-\sigma
\left(  z\right)  }\prod_{k=1}^{K}C\left(  \nu^{k},\underline{\ell}\right)
\rho\left(  z\right)  ^{2\ell-\nu^{k}\cdot\gamma}\\
&  \leq C\left(  \nu,\underline{\ell}\right)  \rho\left(  z\right)
^{2K\ell-\nu\cdot\gamma}e^{-\sigma\left(  z\right)  }\quad,
\end{align*}%
\[
\left\vert \partial^{\nu}\left[  \sigma\left(  z\right)  e^{-\sigma\left(
z\right)  }\right]  \right\vert \leq C\left(  \nu,\underline{\ell}\right)
e^{-\sigma\left(  z\right)  }\cdot\left\{
\begin{array}
[c]{cc}%
\rho\left(  z\right)  ^{2\left(  1+\left\vert \alpha\right\vert \right)
\ell-\nu\cdot\gamma} & ,\quad\text{if }\rho\left(  z\right)  \geq1\ .\\
\rho\left(  z\right)  ^{2\ell-\nu\cdot\gamma} & ,\quad\text{if }\rho\left(
z\right)  \leq1\ .
\end{array}
\right.
\]

Upon combining this bound with (\ref{bound1}) and (\ref{bound2}), while noting
that all multi-indices and linear combinations are determined ultimately by
$\alpha$, $\beta$, and $\mathcal{L}$, we deduce that%
\begin{align}
&  \left\vert \partial^{\alpha}\left[  \varphi_{\varepsilon}\left(  z\right)
\left(  iz\right)  ^{\beta}\sigma(z)e^{-\sigma(z)}L\left(  z\right)
^{-1}\right]  \right\vert \nonumber\\
&  \leq C\left(  \mathcal{L},\alpha,\beta\right)  e^{-\sigma\left(  z\right)
}\cdot\left\{
\begin{array}
[c]{cc}%
\rho\left(  z\right)  ^{2\left\vert \alpha\right\vert \ell+\ell+\beta
\cdot\gamma-\alpha\cdot\gamma} & ,\quad\text{if }\rho\left(  z\right)
\geq1\ ,\\
\rho\left(  z\right)  ^{\ell+\beta\cdot\gamma-\alpha\cdot\gamma} &
,\quad\text{if }0<\rho\left(  z\right)  \leq1\ .
\end{array}
\right. \label{intbound}%
\end{align}
Note that the displayed exponential decay at infinity justifies our previous
integrations by parts. By Lemma \ref{rhoint}, for integrability of this last
expression near $\rho\left(  z\right)  =0$ we require that $\ell+\beta
\cdot\gamma-\alpha\cdot\gamma>-\left\Vert \gamma\right\Vert $. If this
condition holds we may let $\varepsilon\rightarrow0$ inside the integral in
(\ref{jbeta}), as the bounds (\ref{intbound}) are independent of $\varepsilon
$. We have the pointwise limits $\varphi_{\varepsilon}\left(  z\right)
\rightarrow1$ and $\partial^{\nu}\varphi_{\varepsilon}\left(  z\right)
\rightarrow0$ if $\nu\neq0$; thus (\ref{jbeta}) results in%
\begin{equation}
x^{\alpha}J_{\beta}\left(  x,1\right)  =i^{\left\vert \alpha\right\vert }%
\int_{\mathbb{R}^{n}}e^{ix\cdot z}\partial^{\alpha}\left[  \left(  iz\right)
^{\beta}\sigma(z)e^{-\sigma(z)}L\left(  z\right)  ^{-1}\right]  \ dz\quad
,\label{jbetaparts}%
\end{equation}
provided that $\alpha\cdot\gamma<\ell+\beta\cdot\gamma+\left\Vert
\gamma\right\Vert $. Our argument thus far ensures that this integral
converges absolutely, and indeed (\ref{intbound}) and (\ref{jbetaparts}) imply
the bound%
\begin{equation}
\left\vert x^{\alpha}J_{\beta}\left(  x,1\right)  \right\vert \leq C\left(
\mathcal{L},\alpha,\beta\right)  \qquad\text{( if }\alpha\cdot\gamma
<\ell+\beta\cdot\gamma+\left\Vert \gamma\right\Vert \text{ )}\quad
.\label{jay1est}%
\end{equation}

Now in (\ref{jay1est}) we choose $\alpha=N\ell_{k}e_{k}$, where $N$ is a
nonnegative integer and $e_{k}$ is the unit multi-index in the $k$th
coordinate direction. The condition $\alpha\cdot\gamma=N\ell_{k}\gamma
_{k}=N\ell<\ell+\beta\cdot\gamma+\left\Vert \gamma\right\Vert $ leads to the
requirement%
\begin{equation}
0\leq N<1+\frac{\beta\cdot\gamma+\left\Vert \gamma\right\Vert }{\ell}%
\quad.\label{enn}%
\end{equation}
We have $x^{\alpha}=x_{k}{}^{N\ell_{k}}$, and hence (\ref{jay1est}) gives%
\begin{equation}
\left\vert x_{k}{}^{N\ell_{k}}J_{\beta}\left(  x,1\right)  \right\vert \leq
C\left(  \mathcal{L},N,\beta\right)  \quad.\label{jab}%
\end{equation}
If $\delta_{1}$, $\delta_{2}$, $\ldots$, $\delta_{n}$ each have either of the
values $+1$ or $-1$, then application of (\ref{jab}), and (\ref{jay1est}) with
$\alpha=0$, gives%
\[
\left\vert \left(  1+\delta_{1}x_{1}{}^{N\ell_{1}}+\delta_{2}x_{2}{}%
^{N\ell_{2}}+\cdots+\delta_{n}x_{n}{}^{N\ell_{n}}\right)  J_{\beta}\left(
x,1\right)  \right\vert \leq C\left(  \mathcal{L},N,\beta\right)  \quad.
\]
Now, given any $x$ in $\mathbb{R}^{n}$, we may choose the values $\left\{
\delta_{k}\right\}  $ so that this inequality becomes%
\[
\left(  1+\left\vert x_{1}\right\vert ^{N\ell_{1}}+\left\vert x_{2}\right\vert
^{N\ell_{2}}+\cdots+\left\vert x_{n}\right\vert ^{N\ell_{n}}\right)
\left\vert J_{\beta}\left(  x,1\right)  \right\vert \leq C\left(
\mathcal{L},N,\beta\right)  \quad.
\]
But
\[
\left[  1+\rho\left(  x\right)  \right]  ^{N\ell}\leq C\left(  N,\ell
,n\right)  \left(  1+\left\vert x_{1}\right\vert ^{N\ell_{1}}+\left\vert
x_{2}\right\vert ^{N\ell_{2}}+\cdots+\left\vert x_{n}\right\vert ^{N\ell_{n}%
}\right)  \quad,
\]
and thus, provided that (\ref{enn}) holds,
\[
\left\vert J_{\beta}\left(  x,1\right)  \right\vert \leq\frac{C\left(
\mathcal{L},N,\beta\right)  }{\left[  1+\rho\left(  x\right)  \right]
^{N\ell}}\quad\text{.}%
\]
We may choose an integer $N$ satisfying (\ref{enn}) so that $N\geq\left(
\beta\cdot\gamma+\left\Vert \gamma\right\Vert \right)  /\ell$; then we obtain%
\begin{equation}
\left\vert J_{\beta}\left(  x,1\right)  \right\vert \leq\frac{C\left(
\mathcal{L},\beta\right)  }{\left[  1+\rho\left(  x\right)  \right]
^{\beta\cdot\gamma+\left\Vert \gamma\right\Vert }}\quad.\label{jb1}%
\end{equation}

Finally, we set $s=t^{-1/\left(  2\ell\right)  }$ in (\ref{jayz}) to obtain%
\begin{equation}
J_{\beta}\left(  x,t\right)  =t^{-\left(  \beta\cdot\gamma+\left\Vert
\gamma\right\Vert +\ell\right)  /\left(  2\ell\right)  }J_{\beta}\left(
t^{-\gamma/\left(  2\ell\right)  }x,1\right)  \quad.\label{jb2}%
\end{equation}
Applying then (\ref{jb1}), noting that $\rho\left(  t^{-\gamma/\left(
2\ell\right)  }x\right)  =t^{-1/\left(  2\ell\right)  }\rho\left(  x\right)
$, yields (\ref{jaybound}).
\end{proof}

\qquad

Our proposed fundamental solution (\ref{eff}) for $\mathcal{L}$ can be written
in terms of $J$ as the $m\times m$ matrix valued function%
\begin{equation}
F\left(  x\right)  =\left(  2\pi\right)  ^{-n}\int_{0}^{\infty}J\left(
x,t\right)  \ dt\qquad\left(  x\neq0\right)  \quad.\label{effb}%
\end{equation}
In view of (\ref{jayy}), the formal derivative $\partial^{\beta}F$ of
(\ref{effb}) is%
\begin{equation}
F_{\beta}\left(  x\right)  :=\left(  2\pi\right)  ^{-n}\int_{0}^{\infty
}J_{\beta}\left(  x,t\right)  \ dt\qquad\left(  x\neq0\right)  \quad
.\label{effbeta}%
\end{equation}
We will show that, under the added restriction $\left\Vert \gamma\right\Vert
>\ell$, these integrals converge absolutely if $x\neq0$. When $x=0$,
(\ref{effbeta}) and (\ref{jb2}) give%
\[
F_{\beta}\left(  0\right)  =\left(  2\pi\right)  ^{-n}J_{\beta}\left(
0,1\right)  \int_{0}^{\infty}t^{-\left(  \beta\cdot\gamma+\left\Vert
\gamma\right\Vert +\ell\right)  /\left(  2\ell\right)  }\ dt\quad.
\]
As the integral on the right is infinite for any value of $\beta$, formulas
(\ref{effb}) and (\ref{effbeta}) are undefined at $x=0$.

\begin{theorem}
\label{effthm}Suppose $\left\Vert \gamma\right\Vert >\ell$. Then

(a)\ for $x\neq0$ each integral (\ref{effbeta}) for $F_{\beta}\left(
x\right)  $ converges absolutely, and
\begin{equation}
\left\vert F_{\beta}\left(  x\right)  \right\vert \leq C\left(  \mathcal{L}%
,\beta\right)  \rho\left(  x\right)  ^{\ell-\beta\cdot\gamma-\left\Vert
\gamma\right\Vert }\qquad\left(  x\neq0\right)  \ ,\label{effest}%
\end{equation}

(b)\ $F\in C^{\infty}\left(  \mathbb{R}^{n}\backslash\left\{  0\right\}
\right)  $, and $\partial^{\beta}F\left(  x\right)  =F_{\beta}\left(
x\right)  $ for all multi-indices $\beta$ and all nonzero $x$ in
$\mathbb{R}^{n}$,

(c) for $s>0$ and $x\in\mathbb{R}^{n}$,%
\begin{equation}
F_{\beta}\left(  s^{\gamma}x\right)  =s^{\ell-\beta\cdot\gamma-\left\Vert
\gamma\right\Vert }F_{\beta}\left(  x\right)  \quad,\label{effhom}%
\end{equation}

(d)\ $\mathcal{L}F=0$ in the region $\mathbb{R}^{n}\backslash\left\{
0\right\}  $.
\end{theorem}

\begin{proof}
By (\ref{jaybound}),%
\[
\int_{0}^{\infty}\left\vert J_{\beta}\left(  x,t\right)  \right\vert \ dt\leq
C\left(  \mathcal{\mathcal{L}},\beta\right)  \int_{0}^{\infty}\frac{t^{-1/2}%
}{\left[  t^{1/\left(  2\ell\right)  }+\rho\left(  x\right)  \right]
^{\beta\cdot\gamma+\left\Vert \gamma\right\Vert }}\ dt\ .
\]
For $x\neq0$ we make the change of integration parameter $t=\rho\left(
x\right)  ^{2\ell}s^{2\ell}$, to obtain%
\[
\int_{0}^{\infty}\left\vert J_{\beta}\left(  x,t\right)  \right\vert \ dt\leq
C\left(  \mathcal{\mathcal{L}},\beta\right)  \rho\left(  x\right)
^{\ell-\beta\cdot\gamma-\left\Vert \gamma\right\Vert }2\ell\int_{0}^{\infty
}\frac{s^{\ell-1}}{\left(  s+1\right)  ^{\beta\cdot\gamma+\left\Vert
\gamma\right\Vert }}\ ds\quad.
\]
The latter integral converges at zero since $\ell\geq1$, and at infinity since
$\beta\cdot\gamma+\left\Vert \gamma\right\Vert \geq\left\Vert \gamma
\right\Vert >\ell$; thus%
\begin{equation}
\int_{0}^{\infty}\left\vert J_{\beta}\left(  x,t\right)  \right\vert \ dt\leq
C\left(  \mathcal{\mathcal{L}},\beta\right)  \rho\left(  x\right)
^{\ell-\beta\cdot\gamma-\left\Vert \gamma\right\Vert }\qquad\left(
x\neq0\right)  \quad.\label{intjay}%
\end{equation}
Hence (a) follows, with this inequality and (\ref{effbeta}) implying
(\ref{effest}).

To verify the assertions of (b), for $x\neq0$ we examine a difference quotient%
\begin{equation}
\frac{F_{\beta}\left(  x+se_{k}\right)  -F_{\beta}\left(  x\right)  }%
{s}=\left(  2\pi\right)  ^{-n}\int_{0}^{\infty}\frac{J_{\beta}\left(
x+se_{k},t\right)  -J_{\beta}\left(  x,t\right)  }{s}\ dt\quad.\label{dq}%
\end{equation}
If $x\neq0$ and $\left\vert s\right\vert <\left\vert x\right\vert /2$, then
the line connecting $x$ to $x+se_{k}$ misses the origin, and there is a number
$r$ between $0$ and $s$ so that%
\[
\frac{J_{\beta}\left(  x+se_{k},t\right)  -J_{\beta}\left(  x,t\right)  }%
{s}=\frac{\partial}{\partial x_{k}}J_{\beta}\left(  x+re_{k},t\right)
=J_{\beta+e_{k}}\left(  x+re_{k},t\right)  \quad;
\]
then by (\ref{jaybound}),%
\[
\left\vert \frac{J_{\beta}\left(  x+se_{k},t\right)  -J_{\beta}\left(
x,t\right)  }{s}\right\vert \leq C\left(  \mathcal{L},\beta\right)
\frac{t^{-1/2}}{\left[  t^{1/\left(  2\ell\right)  }+\rho\left(
x+re_{k}\right)  \right]  ^{\beta\cdot\gamma+\gamma_{k}+\left\Vert
\gamma\right\Vert }}\ \ .
\]
By the triangle inequality of (\ref{tri}), we have $\rho\left(  x+re_{k}%
\right)  \geq\rho\left(  x\right)  -\rho\left(  re_{k}\right)  \geq\rho\left(
x\right)  /2$ if $s$ (and thus $r$) is sufficiently small, and we obtain%
\[
\left\vert \frac{J_{\beta}\left(  x+se_{k},t\right)  -J_{\beta}\left(
x,t\right)  }{s}\right\vert \leq C\left(  \mathcal{L},\beta\right)
\frac{t^{-1/2}}{\left[  t^{1/\left(  2\ell\right)  }+\rho\left(  x\right)
/2\right]  ^{\beta\cdot\gamma+\gamma_{k}+\left\Vert \gamma\right\Vert }}\ \ .
\]
The condition $\left\Vert \gamma\right\Vert >\ell$ ensures that the right side
of this inequality is an integrable function of $t$ on $\left(  0,\infty
\right)  $; as it is also independent of $s$ we may let $s\rightarrow0$ in
(\ref{dq}) and conclude that%
\[
\frac{\partial}{\partial x_{k}}F_{\beta}\left(  x\right)  =\left(
2\pi\right)  ^{-n}\int_{0}^{\infty}J_{\beta+e_{k}}\left(  x,t\right)
\ dt=F_{\beta+e_{k}}\left(  x\right)  \quad.
\]
An induction arguement now confirms (b).

For $s>0$ and $x\in\mathbb{R}^{n}$, use of (\ref{effbeta}) and (\ref{jayz})
gives%
\[
F_{\beta}\left(  x\right)  =\left(  2\pi\right)  ^{-n}s^{\beta\cdot
\gamma+\left\Vert \gamma\right\Vert +\ell}\int_{0}^{\infty}J_{\beta}\left(
s^{\gamma}x,s^{2\ell}t\right)  \ dt\quad.
\]
In the last integral we make the change of integration parameter $r=s^{2\ell
}t$ to obtain (\ref{effhom}).

To verify statement (d), we use (\ref{op}), the formula $\partial^{\alpha
}F=F_{\alpha}$, (\ref{effbeta}), (\ref{jay}), and (\ref{ell}) to write, for
$x\neq0$,%
\begin{align}
\mathcal{L}F\left(  x\right)   &  =\left(  2\pi\right)  ^{-n}\int_{0}^{\infty
}\int_{\mathbb{R}^{n}}e^{ix\cdot z}L\left(  z\right)  \sigma\left(  z\right)
e^{-t\sigma\left(  z\right)  }L\left(  z\right)  ^{-1}\ dz\ dt\nonumber\\
&  =I\left(  2\pi\right)  ^{-n}\int_{0}^{\infty}\int_{\mathbb{R}^{n}%
}e^{ix\cdot z}\sigma\left(  z\right)  e^{-t\sigma\left(  z\right)
}\ dz\ dt\quad,\label{funda}%
\end{align}
where $I$ is the $m\times m$ identity matrix. It is straightforward to verify
that%
\[
\int_{\mathbb{R}^{n}}e^{ix\cdot z}\sigma\left(  z\right)  e^{-t\sigma\left(
z\right)  }\ dz=-\frac{d}{dt}\int_{\mathbb{R}^{n}}e^{ix\cdot z}e^{-t\sigma
\left(  z\right)  }\ dz\quad,
\]
and so we obtain%
\begin{equation}
\mathcal{L}F\left(  x\right)  =I\left(  2\pi\right)  ^{-n}\left[
\int_{\mathbb{R}^{n}}e^{ix\cdot z}e^{-t\sigma\left(  z\right)  }\ dz\right]
_{t=\infty}^{t=0^{+}}\quad,\label{eval}%
\end{equation}
provided that the evaluations at $t=0^{+}$ and $t=\infty$ exist. To address
this question we consider a scalar valued function%
\[
g_{k}\left(  s,t\right)  =\frac{1}{2\pi}\int_{-\infty}^{\infty}e^{isr}%
e^{-tr^{2k}}\ dr\qquad\left(  s\in\mathbb{R},\text{\ }t>0\right)  \quad,
\]
where $k$ is a positive integer. According to the discussion in chapter 9,
section 2, of the book of Friedman \cite{FR}, $g_{k}\left(  s,t\right)  $ is a
fundamental solution of the parabolic differential equation%
\[
\frac{\partial u\left(  s,t\right)  }{\partial t}=\left(  -1\right)
^{k+1}\frac{\partial^{2k}u\left(  s,t\right)  }{\partial s^{2k}}\quad,
\]
and for $s\in\mathbb{R}$ and $t>0$ satisfies an inequality%
\begin{equation}
\left\vert g_{k}\left(  s,t\right)  \right\vert \leq C_{1}\left(  k\right)
t^{-1/\left(  2k\right)  }\exp\left[  -C_{2}\left(  k\right)  \left(
\frac{s^{2k}}{t}\right)  ^{1/\left(  2k-1\right)  }\right]  \quad,\label{gee}%
\end{equation}
where $C_{1}\left(  k\right)  $ and $C_{2}\left(  k\right)  $ are positive
constants. (See Theorem 1 in chapter 9 of \cite{FR}, or \S 2 of \cite{HAM} for
a more detailed treatment.) In particular, $g_{k}\left(  s,t\right)  $
vanishes at $t=\infty$, and at $t=0^{+}$ provided that $s\neq0$. We now write%
\[
\left(  2\pi\right)  ^{-n}\int_{\mathbb{R}^{n}}e^{ix\cdot z}e^{-t\sigma\left(
z\right)  }\ dz=\left(  2\pi\right)  ^{-n}\int_{\mathbb{R}^{n}}\prod_{j=1}%
^{n}e^{ix_{j}z_{j}}e^{-tz_{j}^{2\ell_{j}}}\ dz=\prod_{j=1}^{n}g_{\ell_{j}%
}\left(  x_{j},t\right)  \ ,
\]%
\[
\left\vert \left(  2\pi\right)  ^{-n}\int_{\mathbb{R}^{n}}e^{ix\cdot
z}e^{-t\sigma\left(  z\right)  }\ dz\right\vert \leq\prod_{j=1}^{n}\left\vert
g_{\ell_{j}}\left(  x_{j},t\right)  \right\vert \quad.
\]
From the bound (\ref{gee}) on $g_{k}$ we deduce that the product on the right
vanishes at $t=\infty$, and at $t=0^{+}$ provided that $x_{j}\neq0$ for some
$j$; thus (\ref{eval}) gives $\mathcal{L}F\left(  x\right)  =0$ if $x\neq0$.
\end{proof}

\qquad

We define an integral operator $\mathcal{S}$, prescribed on suitable
$m\times1$ complex vector functions $f$ on $\mathbb{R}^{n}$ according to%
\begin{equation}
\mathcal{S}f\left(  x\right)  =F\ast f\left(  x\right)  =\int_{\mathbb{R}^{n}%
}F\left(  x-y\right)  \ f\left(  y\right)  \ dy\quad.\label{ess}%
\end{equation}

\begin{theorem}
\label{phithm}Assume $\left\Vert \gamma\right\Vert >\ell$, and let $f$ be an
$m\times1$ complex vector function in the space $C_{0}^{\infty}\left(
\mathbb{R}^{n}\right)  $. Then the integral (\ref{ess}) converges absolutely
for all $x$ in $\mathbb{R}^{n}$, and $\mathcal{S}f\in C^{\infty}\left(
\mathbb{R}^{n}\right)  $ with%
\[
\mathcal{L}\left(  \mathcal{S}f\right)  =f\quad.
\]

\end{theorem}

\begin{proof}
From (\ref{ess}), and (\ref{effest}) with $\beta=0$, we find that%
\[
\left\vert \mathcal{S}f\left(  x\right)  \right\vert \leq C\left(
\mathcal{L}\right)  \int_{\mathbb{R}^{n}}\rho\left(  x-y\right)
^{\ell-\left\Vert \gamma\right\Vert }\ \left\vert f\left(  y\right)
\right\vert \ dy\quad.
\]
Since $f$ has compact support the integral on the right converges at infinity,
and by Lemma \ref{rhoint} it converges near $y=x$ because $\ell-\left\Vert
\gamma\right\Vert >-\left\Vert \gamma\right\Vert $. Thus (\ref{ess}) converges absolutely.

From (\ref{ess}) and (\ref{effb}),%
\begin{equation}
\mathcal{S}f\left(  x\right)  =\left(  2\pi\right)  ^{-n}\int_{\mathbb{R}^{n}%
}\int_{0}^{\infty}J\left(  x-y,t\right)  \ f\left(  y\right)  \ dt\ dy\quad
,\label{essb}%
\end{equation}
and then from (\ref{intjay}) with $\beta=0$,%
\[
\int_{\mathbb{R}^{n}}\int_{0}^{\infty}\left\vert J\left(  x-y,t\right)
\ f\left(  y\right)  \right\vert \ dt\ dy\leq C\left(  \mathcal{\mathcal{L}%
}\right)  \int_{\mathbb{R}^{n}}\rho\left(  x-y\right)  ^{\ell-\left\Vert
\gamma\right\Vert }\ \left\vert f\left(  y\right)  \right\vert \ \ dy\quad.
\]
As again this integral is finite, we may interchange orders of integration in
(\ref{essb}) and substitute (\ref{ja}) to obtain%
\begin{equation}
\mathcal{S}f\left(  x\right)  =\left(  2\pi\right)  ^{-n}\int_{0}^{\infty}%
\int_{\mathbb{R}^{n}}\int_{\mathbb{R}^{n}}e^{i\left(  x-y\right)  \cdot
z}\sigma(z)e^{-t\sigma(z)}L\left(  z\right)  ^{-1}f\left(  y\right)
\ dz\ dy\ dt.\label{essc}%
\end{equation}
Looking at the inner two integrals in (\ref{essc}), we use (\ref{rho}) and
(\ref{rhoineq2}) to estimate%
\begin{align*}
&  \int_{\mathbb{R}^{n}}\int_{\mathbb{R}^{n}}\left\vert e^{i\left(
x-y\right)  \cdot z}\sigma(z)e^{-t\sigma(z)}L\left(  z\right)  ^{-1}f\left(
y\right)  \right\vert \ dz\ dy\\
&  \leq C\left(  \mathcal{L}\right)  \int_{\mathbb{R}^{n}}\int_{\mathbb{R}%
^{n}}\rho(z)^{\ell}\ e^{-t\sigma(z)}\ \left\vert f\left(  y\right)
\right\vert \ dz\ dy\\
&  \leq C\left(  \mathcal{L}\right)  \left\Vert f\right\Vert _{1,\mathbb{R}%
^{n}}\int_{\mathbb{R}^{n}}\rho(z)^{\ell}\ e^{-t\sigma(z)}\ dz<\infty\quad.
\end{align*}
Thus we may interchange orders of integration in these two integrals to write%
\begin{equation}
\mathcal{S}f\left(  x\right)  =\left(  2\pi\right)  ^{-n/2}\int_{0}^{\infty
}\int_{\mathbb{R}^{n}}e^{ix\cdot z}\sigma\left(  z\right)  e^{-t\sigma\left(
z\right)  }L\left(  z\right)  ^{-1}\widehat{f}\left(  z\right)  \ dz\ dt\quad
,\label{essd}%
\end{equation}
where $\widehat{f}$ is the $n$-dimensional Fourier transform of $f$,%
\[
\widehat{f}\left(  z\right)  =\left(  2\pi\right)  ^{-n/2}\int_{\mathbb{R}%
^{n}}e^{-iy\cdot z}\ f\left(  y\right)  \ dy\quad.
\]

Next we use (\ref{rhoineq2}) once more to estimate%
\begin{align*}
&  \int_{0}^{\infty}\int_{\mathbb{R}^{n}}\left\vert e^{ix\cdot z}\sigma\left(
z\right)  e^{-t\sigma\left(  z\right)  }L\left(  z\right)  ^{-1}\widehat
{f}\left(  z\right)  \right\vert \ dz\ dt\\
&  \leq C\left(  \mathcal{L}\right)  \int_{\mathbb{R}^{n}}\rho\left(
z\right)  ^{-\ell}\left\vert \widehat{f}\left(  z\right)  \right\vert \int
_{0}^{\infty}\sigma\left(  z\right)  e^{-t\sigma\left(  z\right)  }\ dt\ dz\\
&  =C\left(  \mathcal{L}\right)  \int_{\mathbb{R}^{n}}\rho\left(  z\right)
^{-\ell}\left\vert \widehat{f}\left(  z\right)  \right\vert \ dz\quad.
\end{align*}
As is well known, if $f\in C_{0}^{\infty}\left(  \mathbb{R}^{n}\right)  $ then
$\left\vert \widehat{f}\left(  z\right)  \right\vert $ decreases at infinity
faster than any power of $\left\vert z\right\vert $. Thus the last integral
converges at infinity, and by Lemma \ref{rhoint} also at zero as we assume
$\ell<\left\Vert \gamma\right\Vert $. Hence we may once more interchange
orders of integration in (\ref{essd}) to arrive at%
\begin{align}
\mathcal{S}f\left(  x\right)   &  =\left(  2\pi\right)  ^{-n/2}\int
_{\mathbb{R}^{n}}e^{ix\cdot z}L\left(  z\right)  ^{-1}\widehat{f}\left(
z\right)  \ \int_{0}^{\infty}\sigma\left(  z\right)  e^{-t\sigma\left(
z\right)  }\ dt\ dz\nonumber\\
&  =\ \left(  2\pi\right)  ^{-n/2}\int_{\mathbb{R}^{n}}e^{ix\cdot z}L\left(
z\right)  ^{-1}\widehat{f}\left(  z\right)  dz\quad.\label{esse}%
\end{align}
It is an easy manner to check that we may differentiate (\ref{esse}) under the
integral to obtain, for any multi-index $\alpha$,%
\begin{equation}
\partial^{\alpha}\mathcal{S}f\left(  x\right)  =\left(  2\pi\right)
^{-n/2}\int_{\mathbb{R}^{n}}e^{ix\cdot z}\left(  iz\right)  ^{\alpha}L\left(
z\right)  ^{-1}\widehat{f}\left(  z\right)  \ dz\quad.\label{essf}%
\end{equation}
Indeed, in view of (\ref{tri}) and (\ref{rhoineq2}), absolute convergence of
these integrals is confirmed by%
\[
\int_{\mathbb{R}^{n}}\left\vert e^{ix\cdot z}\left(  iz\right)  ^{\alpha
}L\left(  z\right)  ^{-1}\widehat{f}\left(  z\right)  \right\vert \ dz\leq
C\left(  \mathcal{L}\right)  \int_{\mathbb{R}^{n}}\rho\left(  z\right)
^{\alpha\cdot\gamma-\ell}\left\vert \widehat{f}\left(  z\right)  \right\vert
\ dz<\infty\quad.
\]

Finally, from (\ref{op}) and (\ref{essf}) it follows that%
\begin{align*}
\mathcal{L}\left(  \mathcal{S}f\right)  \left(  x\right)   &  =\left(
2\pi\right)  ^{-n/2}\int_{\mathbb{R}^{n}}e^{ix\cdot z}L\left(  z\right)
L\left(  z\right)  ^{-1}\widehat{f}\left(  z\right)  \ dz\\
&  =\left(  2\pi\right)  ^{-n/2}\int_{\mathbb{R}^{n}}e^{ix\cdot z}\widehat
{f}\left(  z\right)  \ dz=f\left(  x\right)  \quad,
\end{align*}
with the last equality the Fourier inversion theorem for functions in
$C_{0}^{\infty}\left(  \mathbb{R}^{n}\right)  $.
\end{proof}

\qquad

We mention here related integral representations of Demidenko \cite{DEM2,DEM3}%
, who introduced integral operators $\left\{  P_{h}\right\}  $ defined by%
\begin{align*}
&  \left(  2\pi\right)  ^{n}P_{h}f\left(  x\right) \\
&  =\int_{h}^{h^{-1}}t^{-\left\Vert \gamma\right\Vert /\ell}\int
_{\mathbb{R}^{n}}\int_{\mathbb{R}^{n}}e^{i\left(  x-y\right)  \cdot\left(
t^{-\gamma/\ell}z\right)  }2\kappa\sigma\left(  z\right)  ^{\kappa}%
e^{-\sigma\left(  z\right)  ^{\kappa}}L\left(  z\right)  ^{-1}f\left(
y\right)  \ dz\ dy\ dt\ ,
\end{align*}
where $\kappa$ is a suitable positive integer. Using formulas of
Uspenski\u{\i} \cite{USP} regarding certain averagings of functions, Demidenko
showed that, as $h\rightarrow0$ and under suitable regularity conditions on
$f$, the functions $\left\{  P_{h}f\right\}  $ converge in a weighted Sobolev
norm on $\mathbb{R}^{n}$ to a solution $u$ of $\mathcal{L}u=f$. A modification
of this development leads to the formula for the fundamental solution $F$ and
to the integral operator $\mathcal{S}$ (which, when written as a triple
integral, closely resembles $P_{h}$ after some changes in integration parameters).

\section{Function Spaces}

We introduce function spaces useful in working with semielliptic operators.

Given $0\leq r<\infty$, $1\leq p\leq\infty$, and a domain $\Omega$ in
$\mathbb{R}^{n}$, we say a complex $m\times1$ vector function $u$ is in the
space $W^{r,p}\left(  \Omega,\mathbb{C}^{m},\underline{\ell}\right)  $
provided that $u$ and its weak derivatives $\partial^{\alpha}u$, $0\leq
\alpha\cdot\gamma\leq r$, are in $L^{p}\left(  \Omega\right)  $; the norm of
$u$ in this space is%
\begin{equation}
\left\Vert u\right\Vert _{r,p;\Omega,\underline{\ell}}=\sum_{\alpha\cdot
\gamma\leq r}\left\Vert \partial^{\alpha}u\right\Vert _{p,\Omega}%
\quad.\label{norm1}%
\end{equation}
(We assume always that $\gamma$, $\ell$, and $\underline{\ell}$ are related by
(\ref{gamma}) and (\ref{ellvec}), with $\ell=\max_{k}\ \ell_{k}$.) We say that
$u\in W_{loc}^{r,p}\left(  \Omega,\mathbb{C}^{m},\underline{\ell}\right)  $
whenever $u\in W^{r,p}\left(  \Omega_{0},\mathbb{C}^{m},\underline{\ell
}\right)  $ for all bounded open sets $\Omega_{0}$ with closure in $\Omega$.

If $u$ is defined in all of $\mathbb{R}^{n}$ and $s$ is a real number, we say
$u$ is in the weighted Sobolev space $W_{s}^{r,p}\left(  \mathbb{R}%
^{n},\mathbb{C}^{m},\underline{\ell}\right)  $ provided that the weak
derivatives $\partial^{\alpha}u$, $0\leq\alpha\cdot\gamma\leq r$, are in
$L_{loc}^{p}\left(  \mathbb{R}^{n}\right)  $ and $u$ has finite norm%
\begin{equation}
\left\Vert u\right\Vert _{r,p,s;\underline{\ell}}=\sum_{\alpha\cdot\gamma\leq
r}\left\Vert \left(  1+\rho\right)  ^{s+\alpha\cdot\gamma}\partial^{\alpha
}u\right\Vert _{p,\mathbb{R}^{n}}\quad.\label{norm2}%
\end{equation}
Obviously the spaces $W_{s}^{r,p}\left(  \mathbb{R}^{n},\mathbb{C}%
^{m},\underline{\ell}\right)  $ are decreasing with respect to $s$; that is%
\[
s_{1}\leq s_{2}\Longrightarrow W_{s_{2}}^{r,p}\left(  \mathbb{R}%
^{n},\mathbb{C}^{m},\underline{\ell}\right)  \subset W_{s_{1}}^{r,p}\left(
\mathbb{R}^{n},\mathbb{C}^{m},\underline{\ell}\right)  \quad.
\]
We are concerned in this paper mainly with the cases $r=\ell$ and $r=0$. Note
that in the space $W_{s}^{0,p}\left(  \mathbb{R}^{n},\mathbb{C}^{m}%
,\underline{\ell}\right)  $ the norm (\ref{norm2}) simplifies to%
\[
\left\Vert u\right\Vert _{0,p,s;\underline{\ell}}=\left\Vert \left(
1+\rho\right)  ^{s}u\right\Vert _{p,\mathbb{R}^{n}}\quad.
\]

For $0<R<S\leq\infty$ we define in $\mathbb{R}^{n}$ the bounded open sets%
\begin{equation}
\Omega\left(  R\right)  =\left\{  x:\rho\left(  x\right)  <R\right\}
\qquad,\qquad\Omega\left(  R,S\right)  =\left\{  x:R<\rho\left(  x\right)
<S\right\}  \quad.\label{omega}%
\end{equation}

Following is a density theorem for the spaces $W_{s}^{r,p}\left(
\mathbb{R}^{n},\mathbb{C}^{m},\underline{\ell}\right)  $.

\begin{theorem}
\label{densthm}If\thinspace\ $0\leq r<\infty$, $s\in\mathbb{R}$, and $1\leq
p<\infty$, then $C_{0}^{\infty}\left(  \mathbb{R}^{n}\right)  $ is dense in
the space $W_{s}^{r,p}\left(  \mathbb{R}^{n},\mathbb{C}^{m},\underline{\ell
}\right)  $; that is, given a complex $m\times1$ vector function $u$ in
$W_{s}^{r,p}\left(  \mathbb{R}^{n},\mathbb{C}^{m},\underline{\ell}\right)  $
and $\varepsilon>0$, there exists a complex $m\times1$ vector function
$\varphi$ in $C_{0}^{\infty}\left(  \mathbb{R}^{n}\right)  $ such that%
\begin{equation}
\left\Vert u-\varphi\right\Vert _{r,p,s;\underline{\ell}}=\sum_{\alpha
\cdot\gamma\leq r}\left\Vert \left(  1+\rho\right)  ^{s+\alpha\cdot\gamma
}\partial^{\alpha}\left(  u-\varphi\right)  \right\Vert _{p,\mathbb{R}^{n}%
}<\varepsilon\quad.\label{eps0}%
\end{equation}

\end{theorem}

\begin{proof}
Let $u$ be as described, and suppose $\varepsilon>0$. Let $R$ be a real
constant, $R\geq1$. By Lemma \ref{hmzlemma} there exists a real valued
function $\psi$ in $C_{0}^{\infty}\left(  \mathbb{R}^{n}\right)  $, with
support in the region where $\rho\left(  x\right)  <2R$, such that $0\leq
\psi\leq1$, $\psi\equiv1$ where $\rho\left(  x\right)  \leq R$, and for any
multi-index $\alpha$ and $x$ in $\mathbb{R}^{n}$,%
\begin{equation}
\left\vert \partial^{\alpha}\psi(x)\right\vert \leq C\left(  \ell
,\alpha\right)  R^{-\alpha\cdot\gamma}\quad.\label{psineq}%
\end{equation}
We set $v=\psi u$, so that $v\in W_{s}^{r,p}\left(  \mathbb{R}^{n}%
,\mathbb{C}^{m},\underline{\ell}\right)  $, $v\equiv u$ where $\rho\leq R$,
and $v\equiv0 $ where $\rho\geq2R$. Then for any multi-index $\alpha$ with
$\alpha\cdot\gamma\leq r$,%
\begin{gather}
\left\Vert (1+\rho)^{s+\alpha\cdot\gamma}\partial^{\alpha}(v-u)\right\Vert
_{p,\mathbb{R}^{n}}=\left\Vert (1+\rho)^{s+\alpha\cdot\gamma}\partial^{\alpha
}(v-u)\right\Vert _{p,\Omega(R,\infty)}\nonumber\\
\leq\left\Vert (1+\rho)^{s+\alpha\cdot\gamma}\partial^{\alpha}v\right\Vert
_{p,\Omega(R,2R)}+\left\Vert (1+\rho)^{s+\alpha\cdot\gamma}\partial^{\alpha
}u\right\Vert _{p,\Omega(R,\infty)}\quad.\label{dens1}%
\end{gather}
Use of the product rule for differentiation, along with (\ref{psineq}), gives%
\begin{gather*}
\left\Vert (1+\rho)^{s+\alpha\cdot\gamma}\partial^{\alpha}v\right\Vert
_{p,\Omega(R,2R)}=\left\Vert (1+\rho)^{s+\alpha\cdot\gamma}\sum_{\beta
\leq\alpha}\binom{\alpha}{\beta}\partial^{\beta}u\,\partial^{\alpha-\beta}%
\psi\right\Vert _{p,\Omega(R,2R)}\\
\leq\sum_{\beta\leq\alpha}\left\Vert (1+\rho)^{s+\alpha\cdot\gamma}%
\binom{\alpha}{\beta}\partial^{\beta}u\,C\left(  \ell,\alpha-\beta\right)
R^{-\left(  \alpha-\beta\right)  \cdot\gamma}\right\Vert _{p,\Omega
(R,2R)}\quad.
\end{gather*}
Given the requirement $\alpha\cdot\gamma\leq r$ and $\beta\leq\alpha$, there
are only a finite number of possible values of $\alpha$ and $\beta$ in these
manipulations, depending on $\underline{\ell}$ and $r$. Also, $R\leq1+\rho
\leq3R$ in $\Omega(R,2R)$. It follows that%
\begin{align*}
&  \left\Vert (1+\rho)^{s+\alpha\cdot\gamma}\partial^{\alpha}v\right\Vert
_{p,\Omega(R,2R)}\\
&  \leq C(\underline{\ell},r)\sum_{\beta\leq\alpha}\left\Vert (1+\rho
)^{s+\alpha\cdot\gamma}\partial^{\beta}u\,\left(  1+\rho\right)  ^{-\left(
\alpha-\beta\right)  \cdot\gamma}\right\Vert _{p,\Omega(R,2R)}\\
&  \leq C(\underline{\ell},r)\sum_{\beta\leq\alpha}\left\Vert (1+\rho
)^{s+\beta\cdot\gamma}\,\partial^{\beta}u\right\Vert _{p,\Omega(R,\infty
)}\quad.
\end{align*}
Then from this inequality and (\ref{dens1}),%
\[
\left\Vert \left(  1+\rho\right)  ^{s+\alpha\cdot\gamma}\partial^{\alpha
}(v-u)\right\Vert _{p,\mathbb{R}^{n}}\leq C\left(  \underline{\ell},r\right)
\sum_{\beta\leq\alpha}\left\Vert (1+\rho)^{s+\beta\cdot\gamma}\,\partial
^{\beta}u\right\Vert _{p,\Omega(R,\infty)}\ .
\]
Since the norm (\ref{norm2}) is assumed finite, the right side of this last
inequality tends to $0$ as $R\rightarrow\infty$; thus we may choose $R$ large
enough that%
\begin{equation}
\left\Vert v-u\right\Vert _{r,p,s;\underline{\ell}}=\sum_{\alpha\cdot
\gamma\leq r}\left\Vert \left(  1+\rho\right)  ^{s+\alpha\cdot\gamma}%
\partial^{\alpha}\left(  v-u\right)  \right\Vert _{p,\mathbb{R}^{n}%
}<\varepsilon/2\quad.\label{eps1}%
\end{equation}

Now we use a standard argument involving mollifiers to verify there is a
complex $m\times1$ vector function $\varphi$ in $C_{0}^{\infty}\left(
\mathbb{R}^{n}\right)  $ such that
\begin{equation}
\left\Vert \varphi-v\right\Vert _{r,p,s;\underline{\ell}}<\varepsilon
/2\quad,\label{eps2}%
\end{equation}
which when combined with (\ref{eps1}) yields (\ref{eps0}). Let $\eta$ be a
nonnegative function in $C_{0}^{\infty}\left(  R^{n}\right)  $ vanishing
outside the unit ball $\left\vert x\right\vert \leq1$, with $\int\eta\;dx=1$.
For $t>0$ set $\eta_{t}(x)=t^{-n}\eta(x/t)$, and let $v_{t}$ be the
convolution $v_{t}=\eta_{t}\ast v$. The support of $v$ lies in some ball of
radius $S/2$ centered at $0$, and we may assume $S\geq1$. It follows that
$v_{t}\in C_{0}^{\infty}\left(  R^{n}\right)  $ with support in the ball of
radius $S$ about $0$ if $t<S/2$. For $\alpha\cdot\gamma\leq r$ we have
$\partial^{\alpha}(v_{t})=\left(  \partial^{\alpha}v\right)  _{t}$ and
$\left\Vert \left(  \partial^{\alpha}v\right)  _{t}-\partial^{\alpha
}v\right\Vert _{p,\mathbb{R}^{n}}\rightarrow0$ as $t\rightarrow0$. For
$\left\vert x\right\vert \leq S$ with $S\geq1$, crude estimates yield%
\[
1\leq1+\rho(x)\leq1+nS\quad.
\]
Therefore, for any $\alpha$ with $\alpha\cdot\gamma\leq r$, as $t\rightarrow0
$ we have%
\[
\left\Vert \left(  1+\rho\right)  ^{s+\alpha\cdot\gamma}\partial^{\alpha
}\left(  v_{t}-v\right)  \right\Vert _{p,\mathbb{R}^{n}}\leq C\left(
s,r,S,n\right)  \left\Vert \partial^{\alpha}\left(  v_{t}-v\right)
\right\Vert _{p,\mathbb{R}^{n}}\longrightarrow0\;\;\;.
\]
Thus, if we let $\varphi=v_{t}$ we have (\ref{eps2}) if $t$ is sufficiently small.
\end{proof}

\qquad

Demidenko \cite{DEM1} has also introduced special weighted function spaces for
use with semielliptic operators. He defined a space $W_{p,\tau}^{\underline
{\ell}}\left(  \mathbb{R}^{n}\right)  $, with norm%
\[
\left\Vert u,W_{p,\tau}^{\underline{\ell}}\left(  \mathbb{R}^{n}\right)
\right\Vert =\sum_{\alpha\cdot\gamma\leq\ell}\left\Vert \left(  1+\left\langle
x\right\rangle \right)  ^{-\tau\left(  1-\alpha\cdot\gamma/\ell\right)
}\partial^{\alpha}u\right\Vert _{p,\mathbb{R}^{n}}\quad,
\]
where $\left\langle x\right\rangle $ is defined by%
\[
\left\langle x\right\rangle ^{2}=\sum_{k=1}^{n}x_{k}{}^{2\ell_{k}}=\rho\left(
x\right)  ^{2\ell}\quad.
\]
In terms of $\rho$ this norm is equivalent to%
\[
\sum_{\alpha\cdot\gamma\leq\ell}\left\Vert \left(  1+\rho\right)
^{-\tau\left(  \ell-\alpha\cdot\gamma\right)  }\partial^{\alpha}u\right\Vert
_{p,\mathbb{R}^{n}}\quad.
\]
When $\alpha\cdot\gamma=\ell$ the weight reduces to $1$, regardless of $\tau$;
thus this norm appears fundamentally different from (\ref{norm2}). However, in
the case $\tau=1$, Demidenko's norm corresponds to our norm $\left\Vert
u\right\Vert _{\ell,p,-\ell;\underline{\ell}}$, and the space $W_{p,1}%
^{\underline{\ell}}\left(  \mathbb{R}^{n}\right)  $ is equivalent to
$W_{-\ell}^{\ell,p}\left(  \mathbb{R}^{n},\mathbb{C}^{n},\underline{\ell
}\right)  $. Demidenko has shown \cite{DEM1} that $C_{0}^{\infty}\left(
\mathbb{R}^{n}\right)  $ is dense in $W_{p,\tau}^{\underline{\ell}}\left(
\mathbb{R}^{n}\right)  $ whenever $0\leq\tau\leq1$.

\section{Apriori Bound}

The next result, taken here as a lemma, is a special case of Theorem 2 of
\cite{HMZ}.

\begin{lemma}
\label{apr}Let $\Omega$ be an open subset in $\mathbb{R}^{n}$, and let
$\Omega_{0}$ be a bounded open set whose closure lies in $\Omega$. If
$1<p<\infty$ and $\alpha$ is a multi-index with $\alpha\cdot\gamma\leq\ell$,
then for all complex $m\times1$ functions $u$ in the space $W^{\ell,p}\left(
\Omega,\mathbb{C}^{m},\underline{\ell}\right)  $,%
\[
\left\Vert \partial^{\alpha}u\right\Vert _{p,\Omega_{0}}\leq C\left(
\mathcal{L},p,\Omega,\Omega_{0}\right)  \left[  \left\Vert u\right\Vert
_{p,\Omega}+\left\Vert \mathcal{L}u\right\Vert _{p,\Omega}\right]  \quad.
\]

\end{lemma}

Following is our fundamental apriori bound regarding the operator
$\mathcal{L}$ of (\ref{op}).

\begin{theorem}
\label{apbnd}Let $u$ be a complex $m\times1$ function in the space
$W_{loc}^{\ell,p}\left(  \mathbb{R}^{n},\mathbb{C}^{m},\underline{\ell
}\right)  $. If $1<p<\infty$ and $s\in\mathbb{R}$, then%
\begin{align}
\left\Vert u\right\Vert _{\ell,p,s;\underline{\ell}}  &  =\sum_{\alpha
\cdot\gamma\leq\ell}\left\Vert \left(  1+\rho\right)  ^{s+\alpha\cdot\gamma
}\partial^{\alpha}u\right\Vert _{p,\mathbb{R}^{n}}\nonumber\\
&  \leq C\left(  \mathcal{L},s,p\right)  \left[  \left\Vert \left(
1+\rho\right)  ^{s}u\right\Vert _{p,\mathbb{R}^{n}}+\left\Vert \left(
1+\rho\right)  ^{s+\ell}\mathcal{L}u\right\Vert _{p,\mathbb{R}^{n}}\right]
\quad.\label{apriori}%
\end{align}

\end{theorem}

\begin{proof}
Let $u$, $p$, and $s$ be as described. We assume the right side of
(\ref{apriori}) is finite, as otherwise the inequality is trivial. Let
$\alpha$ be a multi-index such that $\alpha\cdot\gamma\leq\ell$. We use the
notation (\ref{omega}).

First, in the region $\Omega(4)$, where $\rho(x)<4$, Lemma \ref{apr} implies
\[
\left\Vert \partial^{\alpha}u\right\Vert _{p,\Omega(2)}\leq C\left(
\mathcal{L},p\right)  \left[  \left\Vert u\right\Vert _{p,\Omega
(4)}+\left\Vert \mathcal{L}u\right\Vert _{p,\Omega(4)}\right]  \quad.
\]
As $1\leq1+\rho(x)\leq5$ in $\Omega(4)$, this inequality implies%
\begin{align}
&  \int_{\Omega(2)}\left(  1+\rho\right)  ^{\left(  s+\alpha\cdot
\gamma\right)  p}\left\vert \partial^{\alpha}u\right\vert ^{p}\ dx\label{add1}%
\\
&  \leq C\left(  \mathcal{L},s,p\right)  \left[  \int_{\Omega(4)}\left(
1+\rho\right)  ^{sp}\left\vert u\right\vert ^{p}\ dx+\int_{\Omega(4)}\left(
1+\rho\right)  ^{\left(  s+\ell\right)  p}\left\vert \mathcal{L}u\right\vert
^{p}\ dx\right]  \quad.\nonumber
\end{align}

Next let $t$ be a real constant, $t\geq1$, and define a function $v$ by%
\[
v(x)=u\left(  t^{\gamma}x\right)  \quad.
\]
Calculations show that%
\[
\partial^{\alpha}v(x)=t^{\alpha\cdot\gamma}\left(  \partial^{\alpha}u\right)
\left(  t^{\gamma}x\right)  \qquad,\qquad\mathcal{L}v(x)=t^{\ell}\left(
\mathcal{L}u\right)  \left(  t^{\gamma}x\right)  \quad.
\]
Again by Lemma \ref{apr},
\[
\int_{\Omega(2,4)}\left\vert \partial^{\alpha}v(x)\right\vert ^{p}\ dx\leq
C\left(  \mathcal{L},p\right)  \left[  \int_{\Omega(1,8)}\left\vert
v(x)\right\vert ^{p}\ dx+\int_{\Omega(1,8)}\left\vert \mathcal{L}%
v(x)\right\vert ^{p}\ dx\right]  \quad,
\]
or in terms of $u$,%
\begin{align*}
&  \int_{\Omega(2,4)}\left\vert t^{\alpha\cdot\gamma}\left(  \partial^{\alpha
}u\right)  \left(  t^{\gamma}x\right)  \right\vert ^{p}\ dx\\
&  \leq C\left(  \mathcal{L},p\right)  \left[  \int_{\Omega(1,8)}\left\vert
u\left(  t^{\gamma}x\right)  \right\vert ^{p}\ dx+\int_{\Omega(1,8)}\left\vert
t^{\ell}\left(  \mathcal{L}u\right)  \left(  t^{\gamma}x\right)  \right\vert
^{p}\ dx\right]  \quad.
\end{align*}
In these last integrals we make the change of integration parameter
$y=t^{\gamma}x$, with $\rho(y)=t\rho(x)$, $dy=t^{\left\Vert \gamma\right\Vert
}\ dx$, and obtain%
\begin{align*}
&  t^{(\alpha\cdot\gamma)p}\int_{\Omega(2t,4t)}\left\vert \left(
\partial^{\alpha}u\right)  \left(  y\right)  \right\vert ^{p}\ dy\\
&  \leq C\left(  \mathcal{L},p\right)  \left[  \int_{\Omega(t,8t)}\left\vert
u\left(  y\right)  \right\vert ^{p}\ dy+t^{\ell p}\int_{\Omega(t,8t)}%
\left\vert \mathcal{L}u(y)\right\vert ^{p}\ dy\right]  \quad.
\end{align*}
But in $\Omega(t,8t)$ with $t\geq1$, we have $t\leq1+\rho(y)\leq9t$, and so we
may multiply this inequality by $t^{sp}$ to obtain%
\begin{align*}
&  \int_{\Omega(2t,4t)}\left(  1+\rho\right)  ^{\left(  s+\alpha\cdot
\gamma\right)  p}\left\vert \partial^{\alpha}u\right\vert ^{p}\ dy\\
&  \leq C\left(  \mathcal{L},s,p\right)  \left[  \int_{\Omega(t,8t)}\left(
1+\rho\right)  ^{sp}\left\vert u\right\vert ^{p}\ dy+\int_{\Omega
(t,8t)}\left(  1+\rho\right)  ^{\left(  s+\ell\right)  p}\left\vert
\mathcal{L}u\right\vert ^{p}\ dy\right]  \quad.
\end{align*}
Now in this inequality we take $t=2^{m}$ for $m=0,1,2,\ldots$, and add all the
resulting inequalities to (\ref{add1}) to arrive at%
\begin{align*}
&  \int_{\mathbb{R}^{n}}\left(  1+\rho\right)  ^{\left(  s+\alpha\cdot
\gamma\right)  p}\left\vert \partial^{\alpha}u\right\vert ^{p}\ dy\\
&  \leq C\left(  \mathcal{L},s,p\right)  \left[  \int_{\mathbb{R}^{n}}\left(
1+\rho\right)  ^{sp}\left\vert u\right\vert ^{p}\ dy+\int_{\mathbb{R}^{n}%
}\left(  1+\rho\right)  ^{\left(  s+\ell\right)  p}\left\vert \mathcal{L}%
u\right\vert ^{p}\ dy\right]  \quad,
\end{align*}
which leads to%
\[
\left\Vert \left(  1+\rho\right)  ^{s+\alpha\cdot\gamma}\partial^{\alpha
}u\right\Vert _{p,\mathbb{R}^{n}}\leq C\left(  \mathcal{L},s,p\right)  \left[
\left\Vert \left(  1+\rho\right)  ^{s}u\right\Vert _{p,\mathbb{R}^{n}%
}+\left\Vert \left(  1+\rho\right)  ^{s+\ell}\mathcal{L}u\right\Vert
_{p,\mathbb{R}^{n}}\right]  \ .
\]
Finally, we sum over all $\alpha$ such that $\alpha\cdot\gamma\leq\ell$ to
obtain (\ref{apriori}).
\end{proof}

\section{The Operator $\mathcal{S}$}

We investigate properties of the operator $\mathcal{S}$ as a mapping between
certain function spaces. We require another technical lemma.

\begin{lemma}
\label{klemma}For $x$ in $\mathbb{R}^{n}$, and for real numbers $\xi$ and
$\eta$, let%
\[
K\left(  x,\xi,\eta\right)  =\int_{\mathbb{R}^{n}}\rho\left(  x-y\right)
^{\xi}\left[  1+\rho\left(  y\right)  \right]  ^{\eta}\ dy\quad.
\]
If%
\begin{equation}
\xi+\left\Vert \gamma\right\Vert >0\qquad,\qquad\eta+\left\Vert \gamma
\right\Vert >0\qquad,\qquad\xi+\eta+\left\Vert \gamma\right\Vert
<0\quad,\label{conds}%
\end{equation}
then%
\begin{equation}
K\left(  x,\xi,\eta\right)  \leq C\left(  \xi,\eta,\underline{\ell}\right)
\left[  1+\rho\left(  x\right)  \right]  ^{\xi+\eta+\left\Vert \gamma
\right\Vert }\quad.\label{kay}%
\end{equation}

\end{lemma}

\begin{proof}
Note that conditions (\ref{conds}) imply that $\xi,\eta<0$.

Fixing $x$ in $\mathbb{R}^{n}$, we partition $\mathbb{R}^{n}$ into three
disjoint regions,%
\begin{align*}
R_{1}  &  =\left\{  y:\rho\left(  x-y\right)  \leq\frac{1+\rho(x)}{2}\right\}
\quad,\\
R_{2}  &  =\left\{  y:\frac{1+\rho(x)}{2}<\rho\left(  x-y\right)  <2\left[
1+\rho(x)\right]  \right\}  \quad,\\
R_{3}  &  =\left\{  y:2\left[  1+\rho(x)\right]  \leq\rho\left(  x-y\right)
\right\}  \quad,
\end{align*}
and write%
\[
K\left(  x,\xi,\eta\right)  =K_{1}\left(  x,\xi,\eta\right)  +K_{2}\left(
x,\xi,\eta\right)  +K_{3}\left(  x,\xi,\eta\right)  \quad,
\]
where%
\[
K_{i}\left(  x,\xi,\eta\right)  =\int_{R_{i}}\rho\left(  x-y\right)  ^{\xi}
\left[  1+\rho\left(  y\right)  \right]  ^{\eta}\ dy\qquad,\qquad
i=1,2,3\quad.
\]

As $\rho$ satisfies the triangle inequality, in the region $R_{1}$ we have%
\begin{gather*}
1+\rho(x)\leq1+\rho(y)+\rho(x-y)\leq1+\rho(y)+\frac{1+\rho(x)}{2}\quad,\\
1+\rho(x)\leq2\left[  1+\rho(y)\right]  \quad.
\end{gather*}
We use the fact that $\eta<0$, along with Lemma \ref{rhoint}(b) and
$\xi>-\left\Vert \gamma\right\Vert $, to derive%
\begin{align*}
K_{1}(x,\xi,\eta)  &  \leq\left[  \frac{1+\rho\left(  x\right)  }{2}\right]
^{\eta}\int_{\rho\left(  x-y\right)  \leq\left[  1+\rho(x)\right]  /2}%
\rho\left(  x-y\right)  ^{\xi}\ dy\\
&  \leq2^{-\eta}\left[  1+\rho\left(  x\right)  \right]  ^{\eta}\int
_{\rho\left(  z\right)  \leq1+\rho(x)}\rho\left(  z\right)  ^{\xi}\ dz\\
&  =2^{-\eta}\left[  1+\rho\left(  x\right)  \right]  ^{\eta}\left[
1+\rho\left(  x\right)  \right]  ^{\xi+\left\Vert \gamma\right\Vert }%
\int_{\rho\left(  z\right)  \leq1}\rho\left(  z\right)  ^{\xi}\ dz\\
&  =C\left(  \xi,\eta,\underline{\ell}\right)  \left[  1+\rho(x)\right]
^{\xi+\eta+\left\Vert \gamma\right\Vert }\quad.
\end{align*}

In the region $R_{2}$,%
\[
\rho(y)\leq\rho(x-y)+\rho(x)\leq2\left[  1+\rho(x)\right]  +\rho
(x)\leq3\left[  1+\rho(x)\right]  \quad.
\]
We use $\xi<0$ and $-\left\Vert \gamma\right\Vert <\eta<0$, along with Lemma
\ref{rhoint}(b), to derive%
\begin{align*}
K_{2}\left(  x,\xi,\eta\right)   &  \leq\int_{\left[  1+\rho(x)\right]
/2<\rho\left(  x-y\right)  <2\left[  1+\rho(x)\right]  }\left[  \frac
{1+\rho(x)}{2}\right]  ^{\xi}\left[  1+\rho(y)\right]  ^{\eta}\ dy\\
&  \leq2^{-\xi}\left[  1+\rho(x)\right]  ^{\xi}\int_{\rho\left(  y\right)
\leq3\left[  1+\rho(x)\right]  }\rho\left(  y\right)  ^{\eta}\ dy\\
&  =C\left(  \xi,\eta,\underline{\ell}\right)  \left[  1+\rho(x)\right]
^{\xi+\eta+\left\Vert \gamma\right\Vert }\quad.
\end{align*}

In the region $R_{3}$,
\begin{align*}
\rho(x-y)  &  \leq\rho(x)+\rho(y)\leq\frac{\rho(x-y)}{2}+\rho(y)\quad,\\
\rho(x-y)  &  \leq2\rho(y)\leq2\left[  1+\rho(y)\right]  \quad.
\end{align*}
We use the fact that $\eta<0$, along with Lemma \ref{rhoint}(a) and $\xi
+\eta<-\left\Vert \gamma\right\Vert $, to derive%
\begin{align*}
K_{3}\left(  x,\xi,\eta\right)   &  \leq\int_{2\left[  1+\rho(x)\right]
\leq\rho\left(  x-y\right)  }\rho\left(  x-y\right)  ^{\xi}\left[  \frac
{\rho\left(  x-y\right)  }{2}\right]  ^{\eta}\ dy\\
&  \leq2^{-\eta}\int_{1+\rho(x)\leq\rho(z)}\rho(z)^{\xi+\eta}\ dz\leq C\left(
\xi,\eta,\underline{\ell}\right)  \left[  1+\rho(x)\right]  ^{\xi
+\eta+\left\Vert \gamma\right\Vert }\quad.
\end{align*}

Combining finally our estimates for $K_{1}$, $K_{2}$, and $K_{3}$ gives
(\ref{kay}).
\end{proof}

\qquad

\begin{lemma}
\label{esslemma}Suppose $\left\Vert \gamma\right\Vert >\ell$, $1\leq
p\leq\infty$, let $s$ be a real number in the range%
\begin{equation}
\ell-\left\Vert \gamma\right\Vert /p<s<\left\Vert \gamma\right\Vert
-\left\Vert \gamma\right\Vert /p\quad,\label{hyp}%
\end{equation}
and let $f$ be a complex $m\times1$ vector function such that $\left\Vert
\left(  1+\rho\right)  ^{s}f\right\Vert _{p,\mathbb{R}^{n}}<\infty$. Then the
integral%
\begin{equation}
\mathcal{S}f\left(  x\right)  =F\ast f\left(  x\right)  =\int_{\mathbb{R}^{n}%
}F\left(  x-y\right)  \ f\left(  y\right)  \ dy\quad.\label{essit}%
\end{equation}
converges absolutely for almost all $x$ in $\mathbb{R}^{n}$, and%
\begin{equation}
\left\Vert \left(  1+\rho\right)  ^{s-\ell}\mathcal{S}f\right\Vert
_{p,\mathbb{R}^{n}}\leq C\left(  \mathcal{L},p,s\right)  \left\Vert \left(
1+\rho\right)  ^{s}f\right\Vert _{p,\mathbb{R}^{n}}\quad.\label{essbnd}%
\end{equation}
If moreover $p>\left\Vert \gamma\right\Vert /\ell$, then in fact (\ref{essit})
converges absolutely for all $x$ in $\mathbb{R}^{n}$, and%
\begin{equation}
\left\vert \mathcal{S}f\left(  x\right)  \right\vert \leq C\left(
\mathcal{L},p,s\right)  \left[  1+\rho\left(  x\right)  \right]
^{\ell-s-\left\Vert \gamma\right\Vert /p}\left\Vert \left(  1+\rho\right)
^{s}f\right\Vert _{p,\mathbb{R}^{n}}\quad.\label{essbndb}%
\end{equation}

\end{lemma}

\begin{proof}
From (\ref{essit}), and (\ref{effest}) with $\beta=0$,
\begin{align}
\left\vert \mathcal{S}f\left(  x\right)  \right\vert  &  \leq\int
_{\mathbb{R}^{n}}\left\vert F\left(  x-y\right)  \right\vert \ \left\vert
f(y)\right\vert \ dy\nonumber\\
&  \leq C\left(  \mathcal{L}\right)  \int_{\mathbb{R}^{n}}\rho\left(
x-y\right)  ^{\ell-\left\Vert \gamma\right\Vert }\left\vert f\left(  y\right)
\right\vert \ dy,\label{essbd}%
\end{align}

First consider the case $1<p<\infty$. Let $q$ be defined by the usual relation
$1/p+1/q=1$. In general, two finite and nonempty open intervals $\left(
a,b\right)  $ and $\left(  c,d\right)  $ intersect if and only if $a<d$ and
$c<b$. Condition (\ref{hyp}) implies%
\[
\frac{\ell}{q}<s+\frac{\left\Vert \gamma\right\Vert -\ell}{p}\qquad,\qquad
s<\frac{\left\Vert \gamma\right\Vert }{q}\quad;
\]
thus there is a real number $r$ such that%
\begin{equation}
\frac{\ell}{q}<r<\frac{\left\Vert \gamma\right\Vert }{q}\qquad,\qquad
s<r<s+\frac{\left\Vert \gamma\right\Vert -\ell}{p}\quad.\label{intcond}%
\end{equation}
By (\ref{essbd}) and H\"{o}lder's inequality,%
\begin{align*}
\left\vert \mathcal{S}f\left(  x\right)  \right\vert  &  \leq C\left(
\mathcal{L}\right)  \left(  \int_{\mathbb{R}^{n}}\rho\left(  x-y\right)
^{\ell-\left\Vert \gamma\right\Vert }\left[  1+\rho(y)\right]  ^{-rq}%
\ dy\right)  ^{1/q}\\
&  \cdot\left(  \int_{\mathbb{R}^{n}}\rho\left(  x-y\right)  ^{\ell-\left\Vert
\gamma\right\Vert }\left[  1+\rho(y)\right]  ^{rp}\left\vert f(y)\right\vert
^{p}\ dy\right)  ^{1/p}\quad.
\end{align*}
As (\ref{intcond}) and $\ell\geq1$ imply that $\xi=\ell-\left\Vert
\gamma\right\Vert $ and $\eta=-rq$ satisfy the hypotheses of Lemma
\ref{klemma}, we may apply that result and conclude that%
\begin{align*}
\left\vert \mathcal{S}f\left(  x\right)  \right\vert ^{p}  &  \leq C\left(
\mathcal{L},p,s\right)  \left[  1+\rho(x)\right]  ^{-rp+\ell p/q}\\
&  \cdot\int_{\mathbb{R}^{n}}\rho\left(  x-y\right)  ^{\ell-\left\Vert
\gamma\right\Vert }\left[  1+\rho(y)\right]  ^{rp}\left\vert f(y)\right\vert
^{p}\ dy\quad.
\end{align*}
Therefore,%
\begin{align*}
\left[  \left\Vert \left(  1+\rho\right)  ^{s-\ell}\mathcal{S}f\right\Vert
_{p,\mathbb{R}^{n}}\right]  ^{p}  &  =\int_{\mathbb{R}^{n}}\left[
1+\rho\left(  x\right)  \right]  ^{\left(  s-\ell\right)  p}\left\vert
\mathcal{S}f\left(  x\right)  \right\vert ^{p}\ dx\\
&  \leq C\left(  \mathcal{L},p,s\right)  \int_{\mathbb{R}^{n}}\left[
1+\rho(y)\right]  ^{rp}\left\vert f(y)\right\vert ^{p}\\
&  \cdot\int_{\mathbb{R}^{n}}\rho\left(  x-y\right)  ^{\ell-\left\Vert
\gamma\right\Vert }\left[  1+\rho(x)\right]  ^{sp-rp-\ell}\ dx\ dy\ .
\end{align*}
Again, (\ref{intcond}) implies that $\xi=\ell-\left\Vert \gamma\right\Vert $
and $\eta=sp-rp-\ell$ satisfy the hypotheses of Lemma \ref{klemma}; applying
that lemma (with the roles of $x$ and $y$ reversed), we obtain%
\[
\left[  \left\Vert \left(  1+\rho\right)  ^{s-\ell}\mathcal{S}f\right\Vert
_{p,\mathbb{R}^{n}}\right]  ^{p}\leq C\left(  \mathcal{L},p,s\right)
\int_{\mathbb{R}^{n}}\left[  1+\rho(y)\right]  ^{sp}\left\vert f(y)\right\vert
^{p}\ dy\quad,
\]
and thereby (\ref{essbnd}). Note that we have verified that the integral on
the right of (\ref{essbd}) defines a function of $x$ in the space $L_{loc}%
^{p}\left(  \mathbb{R}^{n}\right)  $. In particular, for almost all $x$ this
integral must be finite, and consequently the integral defining $\mathcal{S}%
f\left(  x\right)  $ absolutely convergent.

The case $p=1$, when (\ref{hyp}) reduces to $\ell-\left\Vert \gamma\right\Vert
<s<0$, is simpler. We multiply (\ref{essbd}) by $\left[  1+\rho\left(
x\right)  \right]  ^{s-\ell}$, integrate over $\mathbb{R}^{n}$ with respect to
$x$, and apply Lemma \ref{klemma} as above to obtain (\ref{essbnd}) with $p=1$.

For the case $p=\infty$, when (\ref{hyp}) reduces to $\ell<s<\left\Vert
\gamma\right\Vert $, we apply first (\ref{essbd}) and then Lemma \ref{klemma}
to derive%
\begin{align*}
\left\vert \mathcal{S}f\left(  x\right)  \right\vert  &  \leq C\left(
\mathcal{L}\right)  \left\Vert \left(  1+\rho\right)  ^{s}f\right\Vert
_{\infty}\int_{\mathbb{R}^{n}}\rho\left(  x-y\right)  ^{\ell-\left\Vert
\gamma\right\Vert }\left\vert 1+\rho\left(  y\right)  \right\vert ^{-s}\ dy\\
&  \leq C\left(  \mathcal{L},s\right)  \left\Vert \left(  1+\rho\right)
^{s}f\right\Vert _{\infty}\left[  1+\rho\left(  x\right)  \right]  ^{\ell
-s}\quad.
\end{align*}
This inequality implies (\ref{essbnd}), as well as (\ref{essbndb}), for the
case $p=\infty$.

Next assume $\left\Vert \gamma\right\Vert /\ell<p<\infty$. Application of
H\"{o}lder's inequality to (\ref{essbd}) gives%
\[
\left\vert \mathcal{S}f\left(  x\right)  \right\vert \leq C\left(
\mathcal{L}\right)  \left(  \int_{\mathbb{R}^{n}}\rho\left(  x-y\right)
^{\left(  \ell-\left\Vert \gamma\right\Vert \right)  q}\left[  1+\rho
(y)\right]  ^{-sq}\ dy\right)  ^{1/q}\left\Vert \left(  1+\rho\right)
^{s}f\right\Vert _{p,\mathbb{R}^{n}}\ .
\]
Conditions (\ref{hyp}), and $p>\left\Vert \gamma\right\Vert /\ell$, ensure
that Lemma \ref{klemma} applies with $\xi=\left(  \ell-\left\Vert
\gamma\right\Vert \right)  q$ and $\eta=-sq$, resulting in (\ref{essbndb}).
\end{proof}

\qquad

\begin{theorem}
\label{essthm}Suppose $\left\Vert \gamma\right\Vert >\ell$, $1<p<\infty$, and
let $s$ be a real number in the range%
\begin{equation}
-\left\Vert \gamma\right\Vert /p<s<\left\Vert \gamma\right\Vert -\ell
-\left\Vert \gamma\right\Vert /p\quad.\label{hypb}%
\end{equation}
Let $f$ be a complex $m\times1$ vector function in the space $W_{s+\ell}%
^{0,p}\left(  \mathbb{R}^{n},\mathbb{C}^{m},\underline{\ell}\right)  $; i.e,
such that $\left\Vert \left(  1+\rho\right)  ^{s+\ell}f\right\Vert
_{p,\mathbb{R}^{n}}<\infty$. Then the integral%
\[
\mathcal{S}f\left(  x\right)  =\int_{\mathbb{R}^{n}}F\left(  x-y\right)
\ f\left(  y\right)  \ dy\quad.
\]
converges absolutely for almost all $x$ in $\mathbb{R}^{n}$, and
$\mathcal{S}f\in W_{s}^{\ell,p}\left(  \mathbb{R}^{n},\mathbb{C}%
^{m},\underline{\ell}\right)  $ with $\mathcal{L}\left(  \mathcal{S}f\right)
=f$ and%
\begin{equation}
\sum_{\alpha\cdot\gamma\leq\ell}\left\Vert \left(  1+\rho\right)
^{s+\alpha\cdot\gamma}\partial^{\alpha}\left(  \mathcal{S}f\right)
\right\Vert _{p,\mathbb{R}^{n}}\leq C\left(  \mathcal{L},p,s\right)
\left\Vert \left(  1+\rho\right)  ^{s+\ell}f\right\Vert _{p,\mathbb{R}^{n}%
}\quad.\label{essbndc}%
\end{equation}

\end{theorem}

\begin{proof}
We replace $s$ by $s+\ell$ in Lemma \ref{esslemma}, and deduce that the
integral $\mathcal{S}f\left(  x\right)  $ converges absolutely for almost all
$x$ in $\mathbb{R}^{n}$, with%
\begin{equation}
\left\Vert \left(  1+\rho\right)  ^{s}\mathcal{S}f\right\Vert _{p,\mathbb{R}%
^{n}}\leq C\left(  \mathcal{L},p,s\right)  \left\Vert \left(  1+\rho\right)
^{s+\ell}f\right\Vert _{p,\mathbb{R}^{n}}\quad.\label{essbnd0}%
\end{equation}

By Theorem \ref{densthm} there exists a sequence $\left\{  \varphi
_{k}\right\}  $ of complex $m\times1$ vector functions in $C_{0}^{\infty
}\left(  \mathbb{R}^{n}\right)  $ converging to $f$ in $W_{s+\ell}%
^{0,p}\left(  \mathbb{R}^{n},\mathbb{C}^{m},\underline{\ell}\right)  $, so
that
\begin{equation}
\left\Vert \left(  1+\rho\right)  ^{s+\ell}\left(  \varphi_{k}-f\right)
\right\Vert _{p,\mathbb{R}^{n}}\rightarrow0\quad.\label{phieff}%
\end{equation}
Theorem \ref{phithm} implies $\mathcal{S}\varphi_{k}\in C^{\infty}\left(
\mathbb{R}^{n}\right)  $ for each $k$, with
\begin{equation}
\mathcal{L}\left(  \mathcal{S}\varphi_{k}\right)  =\varphi_{k}\quad
.\label{ellphi}%
\end{equation}
As (\ref{essbnd0}) must apply also to each $\varphi_{k}$ and to $\varphi
_{k}-f$, we have%
\begin{equation}
\left\Vert \left(  1+\rho\right)  ^{s}\mathcal{S}\varphi_{k}\right\Vert
_{p,\mathbb{R}^{n}}\leq C\left(  \mathcal{L},p,s\right)  \left\Vert \left(
1+\rho\right)  ^{s+\ell}\varphi_{k}\right\Vert _{p,\mathbb{R}^{n}}%
\quad,\label{phikay}%
\end{equation}%
\begin{equation}
\left\Vert \left(  1+\rho\right)  ^{s}\mathcal{S}\left(  \varphi_{k}-f\right)
\right\Vert _{p,\mathbb{R}^{n}}\leq C\left(  \mathcal{L},p,s\right)
\left\Vert \left(  1+\rho\right)  ^{s+\ell}\left(  \varphi_{k}-f\right)
\right\Vert _{p,\mathbb{R}^{n}}\quad.\label{phikaf}%
\end{equation}
We apply Theorem \ref{apbnd} to each function $\mathcal{S}\varphi_{k}$ and
obtain%
\begin{align*}
\left\Vert \mathcal{S}\varphi_{k}\right\Vert _{\ell,p,s;\underline{\ell}}  &
=\sum_{\alpha\cdot\gamma\leq\ell}\left\Vert \left(  1+\rho\right)
^{s+\alpha\cdot\gamma}\partial^{\alpha}\left(  \mathcal{S}\varphi_{k}\right)
\right\Vert _{p,\mathbb{R}^{n}}\\
&  \leq C\left(  \mathcal{L},s,p\right)  \left[  \left\Vert \left(
1+\rho\right)  ^{s}\left(  \mathcal{S}\varphi_{k}\right)  \right\Vert
_{p,\mathbb{R}^{n}}+\left\Vert \left(  1+\rho\right)  ^{s+\ell}\varphi
_{k}\right\Vert _{p,\mathbb{R}^{n}}\right]  \quad,
\end{align*}
which when combined with (\ref{phikay}) yields%
\begin{equation}
\sum_{\alpha\cdot\gamma\leq\ell}\left\Vert \left(  1+\rho\right)
^{s+\alpha\cdot\gamma}\partial^{\alpha}\left(  \mathcal{S}\varphi_{k}\right)
\right\Vert _{p,\mathbb{R}^{n}}\leq C\left(  \mathcal{L},s,p\right)
\left\Vert \left(  1+\rho\right)  ^{s+\ell}\varphi_{k}\right\Vert
_{p,\mathbb{R}^{n}}\quad.\label{phikab}%
\end{equation}
But (\ref{phikab}) must apply also to each difference $\varphi_{k}-\varphi
_{j}$, and in view of (\ref{phieff}) we conclude that the sequence $\left\{
\mathcal{S}\varphi_{k}\right\}  $ is Cauchy in the space $W_{s}^{\ell
,p}\left(  \mathbb{R}^{n},\mathbb{C}^{m},\underline{\ell}\right)  $. Hence
$\left\{  \mathcal{S}\varphi_{k}\right\}  $ converges in that space to some
function, which must be $\mathcal{S}f$ because of (\ref{phikaf}) and
(\ref{phieff}). In particular, $\mathcal{S}f\in W_{s}^{\ell,p}\left(
\mathbb{R}^{n},\mathbb{C}^{m},\underline{\ell}\right)  $. Letting
$k\rightarrow\infty$ in (\ref{phikab}) and (\ref{ellphi}) gives (\ref{essbndc}%
) as well as $\mathcal{L}\left(  \mathcal{S}f\right)  =f$.
\end{proof}

\section{Mapping Properties}

We combine our results thus far to draw conclusions about mapping properties
of the partial differential operator
\begin{equation}
\mathcal{L}=\sum_{\alpha\cdot\gamma=\ell}A_{\alpha}\partial^{\alpha}%
\quad,\label{op1}%
\end{equation}
as described in the introduction.

For complex $m\times1$ vector functions $u$ and $v$ on $\mathbb{R}^{n}$, we
define the inner product (when it exists)%
\begin{equation}
\left\langle u,v\right\rangle =\int_{\mathbb{R}^{n}}u\cdot v\ dx=\int
_{\mathbb{R}^{n}}v^{\ast}u\ dx\quad,\label{ip}%
\end{equation}
where \textquotedblleft$^{\ast}$\textquotedblright\ denotes the conjugate
transpose operation.

Recall that functions $u\in W_{s}^{0,p}\left(  \mathbb{R}^{n},\mathbb{C}%
^{m},\underline{\ell}\right)  $ and $v\in W_{-s}^{0,q}\left(  \mathbb{R}%
^{n},\mathbb{C}^{m},\underline{\ell}\right)  $ have the respective norms%
\[
\left\Vert u\right\Vert _{0,p,s;\underline{\ell}}=\left\Vert \left(
1+\rho\right)  ^{s}u\right\Vert _{p,\mathbb{R}^{n}}\qquad,\qquad\left\Vert
v\right\Vert _{0,q,-s;\underline{\ell}}=\left\Vert \left(  1+\rho\right)
^{-s}v\right\Vert _{q,\mathbb{R}^{n}}\quad.
\]
If $1/p+1/q=1$, then (\ref{ip}) is defined for such $u$ and $v$, and according
to H\"{o}lder's inequality,%
\begin{equation}
\left\vert \left\langle u,v\right\rangle \right\vert \leq\left\Vert \left(
1+\rho\right)  ^{s}u\right\Vert _{p,\mathbb{R}^{n}}\left\Vert \left(
1+\rho\right)  ^{-s}v\right\Vert _{q,\mathbb{R}^{n}}=\left\Vert u\right\Vert
_{0,p,s;\underline{\ell}}\left\Vert v\right\Vert _{0,q,-s;\underline{\ell}%
}\quad.\label{hold}%
\end{equation}
Indeed, a standard argument confirms that, if $1<p,q<\infty$ and $1/p+1/q=1$,
then the spaces $W_{s}^{0,p}\left(  \mathbb{R}^{n},\mathbb{C}^{m}%
,\underline{\ell}\right)  $ and $W_{-s}^{0,q}\left(  \mathbb{R}^{n}%
,\mathbb{C}^{m},\underline{\ell}\right)  $ are duals of one another.

If $u\in W_{loc}^{\ell,p}\left(  \mathbb{R}^{n},\mathbb{C}^{m},\underline
{\ell}\right)  $ and $\varphi$ is a complex $m\times1$ vector function in
$C_{0}^{\infty}\left(  \mathbb{R}^{n}\right)  $, we have by the usual
integration by parts that%
\begin{equation}
\left\langle \mathcal{L}u,\varphi\right\rangle =\int_{\mathbb{R}^{n}%
}\mathcal{L}u\cdot\varphi\ dx=\int_{\mathbb{R}^{n}}u\cdot\mathcal{L}^{\ast
}\varphi\ dx=\left\langle u,\mathcal{L}^{\ast}\varphi\right\rangle
\quad,\label{parts1}%
\end{equation}
where $\mathcal{L}^{\ast}$ is the \textit{adjoint operator} to $\mathcal{L}$,%
\[
\mathcal{L}^{\ast}=\sum_{\alpha\cdot\gamma=\ell}\left(  -1\right)
^{\left\vert \alpha\right\vert }A_{\alpha}{}^{\ast}\,\partial^{\alpha}\quad.
\]
Also, as $\left(  \mathcal{L}^{\ast}\right)  ^{\ast}=\mathcal{L}$, we have%
\begin{equation}
\left\langle \mathcal{L}^{\ast}u,\varphi\right\rangle =\left\langle
u,\mathcal{L}\varphi\right\rangle \quad.\label{parts2}%
\end{equation}
If we let $L^{\ast}\left(  x\right)  $ denote the symbol (\ref{ell}) for
$\mathcal{L}^{\ast}$, then by a brief calculation,%
\[
L^{\ast}\left(  x\right)  =L\left(  x\right)  ^{\ast}\quad.
\]
Consequently, $L^{\ast}\left(  x\right)  $ is invertible whenever $L\left(
x\right)  $ is invertible, and semiellipticity of $\mathcal{L}$ implies the
same for $\mathcal{L}^{\ast}$. Moreover, all results proved thus far for
$\mathcal{L}$ are equally valid for $\mathcal{L}^{\ast}$. We let $F^{\ast}$
denote the fundamental solution for the adjoint operator $\mathcal{L}^{\ast}
$,
\[
F^{\ast}\left(  x\right)  =\left(  2\pi\right)  ^{-n}\int_{0}^{\infty}%
\int_{\mathbb{R}^{n}}e^{ix\cdot z}\sigma\left(  z\right)  e^{-t\sigma\left(
z\right)  }L^{\ast}\left(  z\right)  ^{-1}\ dz\ dt\qquad\left(  x\neq0\right)
\quad,
\]
and $\mathcal{S}^{\ast}$ the corresponding convolution operator,%
\[
\mathcal{S}^{\ast}f\left(  x\right)  =F^{\ast}\ast f\left(  x\right)
=\int_{\mathbb{R}^{n}}F^{\ast}\left(  x-y\right)  f\left(  y\right)
\ dy\quad.
\]
Obviously, Theorems \ref{effthm}, \ref{phithm}, and \ref{essthm} apply as well
to $\mathcal{L}^{\ast}$, $F^{\ast}$, and $\mathcal{S}^{\ast}$.

In accordance with (\ref{parts1}) and (\ref{parts2}), for complex $m\times1$
vector functions $u$ and $f$ in $L_{loc}^{1}\left(  \mathbb{R}^{n}\right)  $
we say that $u$ is a \textit{distributional solution} in $\mathbb{R}^{n}$ of
the equation (a)$\ \mathcal{L}u=f$, or (b)$\ \mathcal{L}^{\ast}u=f$, provided
that, respectively,%
\[
\text{(a)\ \ }\left\langle u,\mathcal{L}^{\ast}\varphi\right\rangle
=\left\langle f,\varphi\right\rangle \text{\qquad,\qquad(b)\ \ }\left\langle
u,\mathcal{L}\varphi\right\rangle =\left\langle f,\varphi\right\rangle \quad,
\]
for all complex $m\times1$ vector functions $\varphi$ in $C_{0}^{\infty
}\left(  \mathbb{R}^{n}\right)  $.

Elementary estimates confirm that, for suitable positive constants
$K_{1}\left(  \underline{\ell}\right)  $ and $K_{2}\left(  \underline{\ell
}\right)  $ and for $x\in\mathbb{R}^{n}$,%
\begin{equation}
K_{1}\left(  \underline{\ell}\right)  \left(  1+\left\vert x\right\vert
\right)  ^{1/\ell}\leq1+\rho\left(  x\right)  \leq K_{2}\left(  \underline
{\ell}\right)  \left(  1+\left\vert x\right\vert \right)  \quad.\label{rhoabs}%
\end{equation}
For $z\in\mathbb{C}^{n}$ and $x\in\mathbb{R}^{n}$ we define%
\begin{equation}
\left\Vert z\right\Vert _{x}:=\left[  1+\rho\left(  x\right)  \right]
^{-\ell}\left\vert z\right\vert \quad,\label{temp}%
\end{equation}
and conclude that, for another positive constant $K_{3}\left(  \underline
{\ell}\right)  $,%
\[
K_{3}\left(  \underline{\ell}\right)  \left(  1+\left\vert x\right\vert
\right)  ^{-\ell}\left\vert z\right\vert \leq\left\Vert z\right\Vert _{x}%
\leq\left(  1+\left\vert x\right\vert \right)  ^{\ell}\left\vert z\right\vert
\quad.
\]
This inequality demonstrates, according to the criterion of H\"{o}rmander
(\cite{HO}, \S \ 22.1), that (\ref{temp}) defines a \textit{temperate norm} on
$\mathbb{C}^{n}$, parametrized by $x\in\mathbb{R}^{n}$. For nonzero
$x\in\mathbb{R}^{n}$ and for $z\in\mathbb{C}^{n}$, use of (\ref{rhoineq2})
leads to%
\[
\left\vert z\right\vert =\left\vert L\left(  x\right)  ^{-1}L\left(  x\right)
z\right\vert \leq\left\vert L\left(  x\right)  ^{-1}\right\vert \left\vert
L\left(  x\right)  z\right\vert \leq c_{4}\left(  \mathcal{L}\right)
\rho\left(  x\right)  ^{-\ell}\left\vert L\left(  x\right)  z\right\vert
\quad.
\]
But if $\rho\left(  x\right)  \geq1$ then $\rho\left(  x\right)  ^{-\ell}%
\leq2^{\ell}\left(  1+\rho\left(  x\right)  \right)  ^{-\ell}$, and we obtain%
\begin{equation}
\left\vert z\right\vert \leq C\left(  \mathcal{L}\right)  \left\Vert L\left(
x\right)  z\right\Vert _{x}\qquad\text{,\qquad if }\rho\left(  x\right)
\geq1\quad.\label{crit1}%
\end{equation}
From (\ref{temp}) and (\ref{rhoineq3}) we deduce that%
\begin{align*}
\left\Vert \partial^{\alpha}L\left(  x\right)  z\right\Vert _{x}  &
\leq\left[  1+\rho\left(  x\right)  \right]  ^{-\ell}\left\vert z\right\vert
\cdot\left\{
\begin{tabular}
[c]{cc}%
$c_{5}\left(  \mathcal{L}\right)  \rho\left(  x\right)  ^{\ell-\alpha
\cdot\gamma}$ & ,\ if\ $\alpha\cdot\gamma\leq\ell$,\\
$0$ & ,\ otherwise.
\end{tabular}
\ \ \ \ \right. \\
&  \leq c_{5}\left(  \mathcal{L}\right)  \left[  1+\rho\left(  x\right)
\right]  ^{-\alpha\cdot\gamma}\left\vert z\right\vert \quad.
\end{align*}
Then with use of (\ref{rhoabs}) and the inequality
\[
\alpha\cdot\gamma=\sum_{k=1}^{n}\alpha_{k}\frac{\ell}{\ell_{k}}\geq\sum
_{k=1}^{n}\alpha_{k}=\left\vert \alpha\right\vert \quad,
\]
we find that there is a constant $C\left(  \mathcal{L},\alpha\right)  $ such
that
\begin{equation}
\left\Vert \partial^{\alpha}L\left(  x\right)  z\right\Vert _{x}\leq C\left(
\mathcal{L},\alpha\right)  \left(  1+\left\vert x\right\vert \right)
^{-\left\vert \alpha\right\vert /\ell}\left\vert z\right\vert \quad
,\label{crit2}%
\end{equation}
where obviously $0<1/\ell\leq1$. Inequalities (\ref{crit1}) and (\ref{crit2})
demonstrate that $\mathcal{L}$ is a matrix \textit{hypoelliptic operator}, as
defined by H\"{o}rmander (\cite{HO}, \S \ 22.1). As a consequence (see
\cite{HO}), if $f$ is of class $C^{\infty}$ in an open set in $\mathbb{R}^{n}
$, then any distributional solution of $\mathcal{L}u=f$ in that open set
likewise is of class $C^{\infty}$.

The preceding observations yield the following regularity result.

\begin{proposition}
\label{apdit}Suppose $\left\Vert \gamma\right\Vert >\ell$, $1<p<\infty$,
$s\in\mathbb{R},$ $f\in W_{s+\ell}^{0,p}\left(  \mathbb{R}^{n},\mathbb{C}%
^{m},\underline{\ell}\right)  $, and $u\in W_{s}^{0,p}\left(  \mathbb{R}%
^{n},\mathbb{C}^{m},\underline{\ell}\right)  $. If $u$ is a distributional
solution in $\mathbb{R}^{n}$ of $\mathcal{L}u=f$, then $u\in W_{s}^{\ell
,p}\left(  \mathbb{R}^{n},\mathbb{C}^{m},\underline{\ell}\right)  $, and%
\begin{equation}
\left\Vert u\right\Vert _{\ell,p,s;\underline{\ell}}\leq C\left(
\mathcal{L},s,p\right)  \left[  \left\Vert \left(  1+\rho\right)
^{s}u\right\Vert _{p,\mathbb{R}^{n}}+\left\Vert \left(  1+\rho\right)
^{s+\ell}f\right\Vert _{p,\mathbb{R}^{n}}\right]  \quad.\label{apbd}%
\end{equation}

\end{proposition}

\begin{proof}
Let $B$ be any open ball in $\mathbb{R}^{n}$. As the function $f\chi_{B}$ is
in the space $W_{s^{\prime}+\ell}^{0,p}\left(  \mathbb{R}^{n},\mathbb{C}%
^{m},\underline{\ell}\right)  $ for all $s^{\prime}\in\mathbb{R}$, Theorem
\ref{essthm} asserts that the function $\mathcal{S}\left(  f\chi_{B}\right)  $
is in $W_{s^{\prime}}^{\ell,p}\left(  \mathbb{R}^{n},\mathbb{C}^{m}%
,\underline{\ell}\right)  $ for some $s^{\prime}$, with $\mathcal{L}\left(
\mathcal{S}\left(  f\chi_{B}\right)  \right)  =f\chi_{B}$. Thus
$w=u-\mathcal{S}\left(  f\chi_{B}\right)  $, a distributional solution in $B$
of $\mathcal{L}w=0$, is in $C^{\infty}\left(  B\right)  $. As $B$ is arbitrary
in $\mathbb{R}^{n}$, $u\in W_{loc}^{\ell,p}\left(  \mathbb{R}^{n}%
,\mathbb{C}^{m},\underline{\ell}\right)  $. By Theorem \ref{apbnd}, $u\in
W_{s}^{\ell,p}\left(  \mathbb{R}^{n},\mathbb{C}^{m},\underline{\ell}\right)
$, and (\ref{apbd}) holds.
\end{proof}

\qquad

\begin{lemma}
\label{dual}Assume $\left\Vert \gamma\right\Vert >\ell,$ $1<p<\infty$,
$1/p+1/q=1$, and $s\in\mathbb{R}.$

a)\quad If $u\in W_{s}^{\ell,p}\left(  \mathbb{R}^{n},\mathbb{C}%
^{m},\underline{\ell}\right)  $ and $v\in W_{-s-\ell}^{\ell,q}\left(
\mathbb{R}^{n},\mathbb{C}^{m},\underline{\ell}\right)  $, then%
\begin{equation}
\left\langle \mathcal{L}u,v\right\rangle =\left\langle u,\mathcal{L}^{\ast
}v\right\rangle \qquad,\qquad\left\langle \mathcal{L}^{\ast}u,v\right\rangle
=\left\langle u,\mathcal{L}v\right\rangle \quad.\label{vec1}%
\end{equation}

b)\quad If $f\in W_{-s}^{0,q}\left(  \mathbb{R}^{n},\mathbb{C}^{m}%
,\underline{\ell}\right)  $ and
\begin{equation}
-\frac{\left\Vert \gamma\right\Vert }{p}<s<\frac{\left\Vert \gamma\right\Vert
}{q}-\ell\quad,\label{conb}%
\end{equation}
then $\mathcal{S}f$ and $\mathcal{S}^{\ast}f$ are in $W_{-s-\ell}^{\ell
,q}\left(  \mathbb{R}^{n},\mathbb{C}^{m},\underline{\ell}\right)  $, with
$\mathcal{L}\left(  \mathcal{S}f\right)  =\mathcal{L}^{\ast}\left(
\mathcal{S}^{\ast}f\right)  =f$; moreover, for $u\in W_{s}^{\ell,p}\left(
\mathbb{R}^{n},\mathbb{C}^{m},\underline{\ell}\right)  $,%
\begin{equation}
\left\langle \mathcal{L}u,\mathcal{S}^{\ast}f\right\rangle =\left\langle
\mathcal{L}^{\ast}u,\mathcal{S}f\right\rangle =\left\langle u,f\right\rangle
\quad.\label{vec2}%
\end{equation}

\end{lemma}

\begin{proof}
a)\quad By Theorem \ref{densthm}, there exists a sequence $\left\{
\varphi_{k}\right\}  $ of complex $m\times1$ vector functions in
$C_{0}^{\infty}\left(  \mathbb{R}^{n}\right)  $ converging to $v$ in the space
$W_{-s-\ell}^{\ell,q}\left(  \mathbb{R}^{n},\mathbb{C}^{m},\underline{\ell
}\right)  $; that is, with%
\[
\sum_{\alpha\cdot\gamma\leq\ell}\left\Vert \left(  1+\rho\right)
^{-s-\ell+\alpha\cdot\gamma}\partial^{\alpha}\left(  v-\varphi_{k}\right)
\right\Vert _{q,\mathbb{R}^{n}}\longrightarrow0\quad.
\]
Since $u\in W_{loc}^{\ell,p}\left(  \mathbb{R}^{n},\mathbb{C}^{m}%
,\underline{\ell}\right)  $, for each $\varphi_{k}$ we have%
\[
\left\langle \mathcal{L}u,\varphi_{k}\right\rangle =\left\langle
u,\mathcal{L}^{\ast}\varphi_{k}\right\rangle \quad.
\]
Since also $\left\Vert \left(  1+\rho\right)  ^{s}u\right\Vert _{p,\mathbb{R}%
^{n}}<\infty$ and $\left\Vert \left(  1+\rho\right)  ^{s+\ell}\mathcal{L}%
u\right\Vert _{p,\mathbb{R}^{n}}<\infty$, in view of (\ref{hold}) we may let
$k\rightarrow\infty$ in this equation to obtain the left equation of
(\ref{vec1}). The right equation of (\ref{vec1}) follows similarly, or by
replacing $\mathcal{L}$ with $\mathcal{L}^{\ast}$.

b)\quad We apply Theorem \ref{essthm} but with $p$ replaced by $q$ and $s$
replaced by $-s-\ell$. The hypothesis (\ref{hypb}) is replaced by
\[
-\left\Vert \gamma\right\Vert /q<-s-\ell<\left\Vert \gamma\right\Vert
-\ell-\left\Vert \gamma\right\Vert /q\quad,
\]
which follows from (\ref{conb}). The theorem concludes that $\mathcal{S}f$ and
$\mathcal{S}^{\ast}f$ are in $W_{-s-\ell}^{\ell,q}\left(  \mathbb{R}%
^{n},\mathbb{C}^{m},\underline{\ell}\right)  $, with $\mathcal{L}\left(
\mathcal{S}f\right)  =\mathcal{L}^{\ast}\left(  \mathcal{S}^{\ast}f\right)
=f$. Application of (\ref{vec1}) to $v=\mathcal{S}^{\ast}f$ and $v=\mathcal{S}%
f$ yields (\ref{vec2}).
\end{proof}

\begin{theorem}
\label{main}Assume $\left\Vert \gamma\right\Vert >\ell$, $1<p<\infty$,
$1/p+1/q=1$, $s\in\mathbb{R}$, and consider the mapping%
\begin{equation}
\mathcal{L}:W_{s}^{\ell,p}\left(  \mathbb{R}^{n},\mathbb{C}^{m},\underline
{\ell}\right)  \longrightarrow W_{s+\ell}^{0,p}\left(  \mathbb{R}%
^{n},\mathbb{C}^{m},\underline{\ell}\right)  \quad.\label{map}%
\end{equation}

a)\quad If $-\left\Vert \gamma\right\Vert /p<s$, then the mapping is one-to-one.

b)\quad If $-\left\Vert \gamma\right\Vert /p<s<\left\Vert \gamma\right\Vert
/q-\ell$, the mapping is onto,

c)\quad If $s<-\left\Vert \gamma\right\Vert /p$, the mapping is not one-to-one.

d)\quad If $s\geq\left\Vert \gamma\right\Vert /q-\ell$ the mapping is not onto.

e)\quad If $s=-\left\Vert \gamma\right\Vert /p$, the mapping is not bounded below.

Consequently, (\ref{map}) is an isomorphism if and only if
\begin{equation}
-\left\Vert \gamma\right\Vert /p<s<\left\Vert \gamma\right\Vert /q-\ell
\quad.\label{iso}%
\end{equation}

\end{theorem}

\begin{proof}
a)\quad Under the given assumptions, assume $u\in W_{s}^{\ell,p}\left(
\mathbb{R}^{n},\mathbb{C}^{m},\underline{\ell}\right)  $ and $\mathcal{L}u=0$.
We may choose $s$ smaller, if necessary, so that (\ref{conb}) holds. Then by
Lemma \ref{dual}(b), for all $f$ in $W_{-s}^{0,q}\left(  \mathbb{R}%
^{n},\mathbb{C}^{m},\underline{\ell}\right)  $ we have%
\[
\left\langle u,f\right\rangle =\left\langle \mathcal{L}u,\mathcal{S}^{\ast
}f\right\rangle =\left\langle 0,\mathcal{S}^{\ast}f\right\rangle =0\quad.
\]
As $u\in W_{s}^{0,p}\left(  \mathbb{R}^{n},\mathbb{C}^{m},\underline{\ell
}\right)  $ and $W_{-s}^{0,q}\left(  \mathbb{R}^{n},\mathbb{C}^{m}%
,\underline{\ell}\right)  $ is the dual space, we infer that $u=0$. Thus
(\ref{map}) is one-to-one.

b)\quad Under the given assumptions, let $f\in W_{s+\ell}^{0,p}\left(
\mathbb{R}^{n},\mathbb{C}^{m},\underline{\ell}\right)  $. By Theorem
\ref{essthm}, $\mathcal{S}f\in W_{s}^{\ell,p}\left(  \mathbb{R}^{n}%
,\mathbb{C}^{m},\underline{\ell}\right)  $ and $\mathcal{L}\left(
\mathcal{S}f\right)  =f$. Thus $\mathcal{L}$ is onto.

c)\quad Since $\mathcal{L}$ has no zero order term, any constant vector $u$
solves $\mathcal{L}u=0$. But for such $u$, (\ref{norm2}) gives%
\[
\left\Vert u\right\Vert _{\ell,p,s;\underline{\ell}}=\left\Vert \left(
1+\rho\right)  ^{s}u\right\Vert _{p,\mathbb{R}^{n}}=\left\vert u\right\vert
\left(  \int_{\mathbb{R}^{n}}\left(  1+\rho\right)  ^{sp}\ dx\right)
^{1/p}\quad,
\]
which according to Lemma \ref{rhoint} is finite whenever $sp<-\left\Vert
\gamma\right\Vert $. Thus, if $s<-\left\Vert \gamma\right\Vert /p$, then $u\in
W_{s}^{\ell,p}\left(  \mathbb{R}^{n},\mathbb{C}^{m},\underline{\ell}\right)  $
and (\ref{map}) is not one-to-one.

d)\quad We consider the case $s>\left\Vert \gamma\right\Vert /q-\ell$, but
delay the case $s=\left\Vert \gamma\right\Vert /q-\ell$ until after the proof
of (e). Let $v$ be any nonzero $m\times1$ constant function, and let $f $ be
the function%
\[
f\left(  x\right)  =\left[  1+\rho\left(  x\right)  \right]  ^{-\left(
s+\ell+\left\Vert \gamma\right\Vert \right)  }v\quad.
\]
By Lemma \ref{rhoint},%
\[
\left\Vert v\right\Vert _{\ell,q,-s-\ell;\underline{\ell}}=\left\vert
v\right\vert \left(  \int_{\mathbb{R}^{n}}\left(  1+\rho\right)  ^{-\left(
s+\ell\right)  q}\ dx\right)  ^{1/q}<\infty\quad,
\]
since $-\left(  s+\ell\right)  q<-\left\Vert \gamma\right\Vert $; also,%
\begin{align*}
\left\Vert f\right\Vert _{0,p,s+\ell;\underline{\ell}}  &  =\left(
\int_{\mathbb{R}^{n}}\left(  1+\rho\right)  ^{\left(  s+\ell\right)
p}\left\vert f\right\vert ^{p}\ dx\right)  ^{1/p}\\
&  =\left\vert v\right\vert \left(  \int_{\mathbb{R}^{n}}\left(
1+\rho\right)  ^{-\left\Vert \gamma\right\Vert p}\ dx\right)  ^{1/p}%
<\infty\quad,
\end{align*}
as $-\left\Vert \gamma\right\Vert p<-\left\Vert \gamma\right\Vert $. Thus
$v\in W_{-s-\ell}^{\ell,q}\left(  \mathbb{R}^{n},\mathbb{C}^{m},\underline
{\ell}\right)  $ and $f\in W_{s+\ell}^{0,p}\left(  \mathbb{R}^{n}%
,\mathbb{C}^{m},\underline{\ell}\right)  $. If $f=\mathcal{L}u$ for some $u$
in $W_{s}^{\ell,p}\left(  \mathbb{R}^{n},\mathbb{C}^{m},\underline{\ell
}\right)  $, then Lemma \ref{dual}(a) requires that
\begin{align*}
0  &  =\left\langle u,0\right\rangle =\left\langle u,\mathcal{L}^{\ast
}v\right\rangle =\left\langle \mathcal{L}u,v\right\rangle =\left\langle
f,v\right\rangle \\
&  =\int_{\mathbb{R}^{n}}f\cdot v\ dx=\left\vert v\right\vert ^{2}%
\int_{\mathbb{R}^{n}}\left[  1+\rho\left(  x\right)  \right]  ^{-\left(
s+\ell+\left\Vert \gamma\right\Vert \right)  }\ dx>0\quad,
\end{align*}
a contradiction. Thus (\ref{map}) is not onto.

e)\quad Let $s=-\left\Vert \gamma\right\Vert /p$. Given $R\geq1$, let $\psi$
be the function described in the proof of Theorem \ref{densthm}, satisfying
$\psi\in C_{0}^{\infty}\left(  \mathbb{R}^{n}\right)  $, $0\leq\psi\leq1$,
$\psi\left(  x\right)  =1$ for $\rho\left(  x\right)  \leq R$, $\psi\left(
x\right)  =0$ for $\rho\left(  x\right)  \geq2R$, and
\[
\left\vert \partial^{\alpha}\psi\left(  x\right)  \right\vert \leq C\left(
\ell,\alpha\right)  R^{-\alpha\cdot\gamma}\quad.
\]
Let $v$ be any nonzero constant $m\times1$ vector, and let $u$ be the function
$u\left(  x\right)  =\psi\left(  x\right)  v$. Then%
\[
\left\vert \mathcal{L}u\left(  x\right)  \right\vert =\left\vert \sum
_{\alpha\cdot\gamma=\ell}A_{\alpha}\partial^{\alpha}\psi\left(  x\right)
v\right\vert \leq\left\{
\begin{tabular}
[c]{ll}%
$C\left(  \mathcal{L}\right)  \left\vert v\right\vert R^{-\ell}$ & , if
$R\leq\rho\left(  x\right)  \leq2R\,,$\\
$0$ & , otherwise\thinspace.
\end{tabular}
\ \ \ \right.
\]
As $1\leq R\leq\rho\leq2R$ implies $\rho\leq1+\rho\leq2\rho$, we have%
\[
\left\Vert \left(  1+\rho\right)  ^{s+\ell}\mathcal{L}u\right\Vert
_{p,\mathbb{R}^{n}}\leq C\left(  \mathcal{L}\right)  \left\vert v\right\vert
R^{-\ell}\left(  \int_{R\leq\rho\left(  x\right)  \leq2R}\rho\left(  x\right)
^{\left(  s+\ell\right)  p}\ dx\right)  ^{1/p}\quad.
\]
In the last integral we set $x=R^{\gamma}y$, $\rho\left(  x\right)
=R\rho\left(  y\right)  $, $dx=R^{\left\Vert \gamma\right\Vert }dy$ and obtain%
\begin{align*}
\int_{R\leq\rho\left(  x\right)  \leq2R}\rho\left(  x\right)  ^{\left(
s+\ell\right)  p}\ dx  &  =R^{\left(  s+\ell\right)  p+\left\Vert
\gamma\right\Vert }\int_{1\leq\rho\left(  y\right)  \leq2}\rho\left(
y\right)  ^{\left(  s+\ell\right)  p}\ dy\\
&  =C\left(  \mathcal{L},p\right)  R^{\left(  s+\ell\right)  p+\left\Vert
\gamma\right\Vert }\quad,
\end{align*}
and thereby, upon setting $s=-\left\Vert \gamma\right\Vert /p$,%
\begin{equation}
\left\Vert \left(  1+\rho\right)  ^{s+\ell}\mathcal{L}u\right\Vert
_{p,\mathbb{R}^{n}}\leq C\left(  \mathcal{L},p\right)  \left\vert v\right\vert
\quad.\label{ellu}%
\end{equation}
On the other hand, Lemma \ref{rhoint} ensures that, as $R\rightarrow\infty$,%
\begin{align*}
\left\Vert \left(  1+\rho\right)  ^{s}u\right\Vert _{p,\mathbb{R}^{n}}  &
\geq\left(  \int_{\rho\left(  x\right)  \leq R}\left[  1+\rho\left(  x\right)
\right]  ^{sp}\left\vert u\left(  x\right)  \right\vert ^{p}\ dx\right)
^{1/p}\\
&  =\left\vert v\right\vert \left(  \int_{\rho\left(  x\right)  \leq R}\left[
1+\rho\left(  x\right)  \right]  ^{-\left\Vert \gamma\right\Vert }\ dx\right)
^{1/p}\longrightarrow\infty\quad.
\end{align*}
Combining this result with (\ref{ellu}), we conclude that, as $R\rightarrow
\infty$,%
\[
\frac{\left\Vert \mathcal{L}u\right\Vert _{0,p,s+\ell;\underline{\ell}}%
}{\left\Vert u\right\Vert _{\ell,p,s;\underline{\ell}}}\leq\frac{C\left(
\mathcal{L},p\right)  \left\vert v\right\vert }{\left\Vert \left(
1+\rho\right)  ^{s}u\right\Vert _{p,\mathbb{R}^{n}}}\longrightarrow0\quad.
\]
By choosing $R$ large enough we can make the left side of this expression as
small as desired; thus the mapping (\ref{map}) is not bounded below in the
case $s=-\left\Vert \gamma\right\Vert /p$.

We are left with only the case $s=\left\Vert \gamma\right\Vert /q-\ell$. We
assume (\ref{map}) is onto, and we will produce a contradiction. The condition
$\left\Vert \gamma\right\Vert >\ell$ implies $s>-\left\Vert \gamma\right\Vert
/p$, and the result of (a) confirms the mapping is one-to-one. Consequently,
as a bounded, onto, one-to-one.mapping from one Banach space to another,
$\mathcal{L}$ has a bounded inverse mapping $\mathcal{L}^{-1}$. In
particular,
\[
\mathcal{L}^{-1}:W_{s+\ell}^{0,p}\left(  \mathbb{R}^{n},\mathbb{C}%
^{m},\underline{\ell}\right)  \longrightarrow W_{s}^{\ell,p}\left(
\mathbb{R}^{n},\mathbb{C}^{m},\underline{\ell}\right)
\]
and there is a positive constant $M$ such that%
\[
\left\Vert \mathcal{L}^{-1}v\right\Vert _{\ell,p,s;\underline{\ell}}\leq
M\left\Vert v\right\Vert _{0,p,s+\ell;\underline{\ell}}\quad.
\]

Given any function $f$ in $W_{-s}^{0,q}\left(  \mathbb{R}^{n},\mathbb{C}%
^{m},\underline{\ell}\right)  $, we may define a mapping $T$ from $W_{s+\ell
}^{0,p}\left(  \mathbb{R}^{n},\mathbb{C}^{m},\underline{\ell}\right)  $ into
$\mathbb{C}$, according to%
\[
T\left(  v\right)  =\left\langle \mathcal{L}^{-1}v,f\right\rangle \quad.
\]
Obviously $T$ is linear, and it is also bounded, as verified by%
\begin{align*}
\left\vert T\left(  v\right)  \right\vert  &  =\left\vert \left\langle
\mathcal{L}^{-1}v,f\right\rangle \right\vert \leq\left\Vert \mathcal{L}%
^{-1}v\right\Vert _{0,p,s;\underline{\ell}}\left\Vert f\right\Vert
_{0,q,-s;\underline{\ell}}\\
&  \leq\left\Vert \mathcal{L}^{-1}v\right\Vert _{\ell,p,s;\underline{\ell}%
}\left\Vert f\right\Vert _{0,q,-s;\underline{\ell}}\leq M\left\Vert
f\right\Vert _{0,q,-s;\underline{\ell}}\left\Vert v\right\Vert _{0,p,s+\ell
;\underline{\ell}}\quad.
\end{align*}
Since $W_{-s-\ell}^{0,q}\left(  \mathbb{R}^{n},\mathbb{C}^{m},\underline{\ell
}\right)  $ is the dual of $W_{s+\ell}^{0,p}\left(  \mathbb{R}^{n}%
,\mathbb{C}^{m},\underline{\ell}\right)  $, there exists a function $u$ in
$W_{-s-\ell}^{0,q}\left(  \mathbb{R}^{n},\mathbb{C}^{m},\underline{\ell
}\right)  $ such that%
\[
T\left(  v\right)  =\left\langle \mathcal{L}^{-1}v,f\right\rangle
=\left\langle v,u\right\rangle \qquad,\qquad\forall v\in W_{s+\ell}%
^{0,p}\left(  \mathbb{R}^{n},\mathbb{C}^{m},\underline{\ell}\right)  \quad.
\]
Giving any complex $m\times1$ vector function $\varphi$ in $C_{0}^{\infty
}\left(  \mathbb{R}^{n}\right)  $ we may choose $v=\mathcal{L}\varphi$ to
obtain%
\[
\left\langle \varphi,f\right\rangle =\left\langle \mathcal{L}\varphi
,u\right\rangle \qquad,\qquad\left\langle f,\varphi\right\rangle =\left\langle
u,\mathcal{L}\varphi\right\rangle \quad.
\]
This relation implies that $u$ is a distributional solution in $\mathbb{R}%
^{n}$ of the equation $\mathcal{L}^{\ast}u=f$. By Proposition \ref{apdit}
applied to $\mathcal{L}^{\ast}$, and with $p$ replaced by $q$ and $s$ by
$-s-\ell$, we have $u\in W_{-s-\ell}^{\ell,q}\left(  \mathbb{R}^{n}%
,\mathbb{C}^{m},\underline{\ell}\right)  $. As $f$ is arbitrary in
$W_{-s}^{0,q}\left(  \mathbb{R}^{n},\mathbb{C}^{m},\underline{\ell}\right)  $,
this argument shows that the mapping
\[
\mathcal{L}^{\ast}:W_{-s-\ell}^{\ell,q}\left(  \mathbb{R}^{n},\mathbb{C}%
^{m},\underline{\ell}\right)  \longrightarrow W_{-s}^{0,q}\left(
\mathbb{R}^{n},\mathbb{C}^{m},\underline{\ell}\right)
\]
is onto. But the result of (e), with $\mathcal{L}$ replaced by $\mathcal{L}%
^{\ast}$, $p$ by $q$, and $s$ by $-s-\ell$, applies to this mapping, as
$-s-\ell=-\left\Vert \gamma\right\Vert /q$. Since $\mathcal{L}^{\ast}$ is onto
but not bounded below, it cannot be one-to-one. Hence there is a function $w$
in $W_{-s-\ell}^{\ell,q}\left(  \mathbb{R}^{n},\mathbb{C}^{m},\underline{\ell
}\right)  $ such that $w\neq0$ and $\mathcal{L}^{\ast}w=0$. By Lemma
\ref{dual}(a) we then have, for all $u$ in $W_{s}^{\ell,p}\left(
\mathbb{R}^{n},\mathbb{C}^{m},\underline{\ell}\right)  $,%
\[
\left\langle \mathcal{L}u,w\right\rangle =\left\langle u,\mathcal{L}^{\ast
}w\right\rangle =\left\langle u,0\right\rangle =0\quad.
\]
But we assume (\ref{map}) is onto, so we have $\left\langle f,w\right\rangle
=0$ for all $f$ in $W_{s+\ell}^{0,p}\left(  \mathbb{R}^{n},\mathbb{C}%
^{m},\underline{\ell}\right)  $, a contradiction if $w\neq0$.
\end{proof}

\qquad

As a special case of Theorem \ref{main}, we take $s=-\ell$ and find that the
mapping%
\begin{equation}
\mathcal{L}:W_{-\ell}^{\ell,p}\left(  \mathbb{R}^{n}\mathbb{C}^{m}%
,\underline{\ell}\right)  \longrightarrow W_{0}^{0,p}\left(  \mathbb{R}%
^{n},\mathbb{C}^{m},\underline{\ell}\right) \label{dimi}%
\end{equation}
is an isomorphism provided that $\ell<\left\Vert \gamma\right\Vert /p$. This
result has already been obtained by Demidenko \cite{DEM2,DEM3}, who wrote
(\ref{dimi}) with the notation%
\[
\mathcal{L}:W_{p,1}^{\underline{\ell}}\left(  \mathbb{R}^{n}\right)
\longrightarrow L_{p}\left(  \mathbb{R}^{n}\right)  \quad.
\]

\section{Examples}

We give a few examples to which results of the paper apply.

\begin{example}
Consider a parabolic opertor in $\mathbb{R}^{n+1}=\mathbb{R}^{n}%
\times\mathbb{R}$,%
\begin{equation}
\mathcal{L}u=\sum_{\left\vert \alpha\right\vert =\ell}A_{\alpha}\partial
_{x}^{\alpha}u-I\partial_{t}u\quad,\label{par}%
\end{equation}
where each $A_{\alpha}$ is a complex constant $m\times m$ matrix, $I$ is the
$m\times m$ identity, and $u=u\left(  x_{1},\ldots,x_{n},t\right)  $. The
usual parabolicity condition (see \cite{FR}) requires that each eigenvalue
$\lambda\left(  x\right)  $ of the matrix%
\[
P\left(  x\right)  =\sum_{\left\vert \alpha\right\vert =\ell}A_{\alpha}\left(
ix\right)  ^{\alpha}%
\]
satisfy an inequality%
\begin{equation}
\operatorname{Re}\ \lambda\left(  x\right)  \leq-\delta\left\vert x\right\vert
^{\ell}\qquad\left(  \delta>0\right)  \quad.\label{par1}%
\end{equation}
For (\ref{par}), formula (\ref{ell}) gives%
\[
L\left(  x,t\right)  =\sum_{\left\vert \alpha\right\vert =\ell}A_{\alpha
}\left(  ix\right)  ^{\alpha}-itI\quad.
\]
Thus $L\left(  x,t\right)  $ is invertible if $\left(  x,t\right)  \neq0$, as
(\ref{par1}) shows that $P\left(  x\right)  $ has no purely imaginary
eigenvalue when $x\neq0$. Also for (\ref{par}), whose order is $\ell$, we
determine that%
\[
\underline{\ell}=\left(  \ell,\ldots,\ell,1\right)  \qquad,\qquad
\gamma=\left(  1,\ldots,1,\ell\right)  \qquad,\qquad\left\Vert \gamma
\right\Vert =n+\ell>\ell\quad,
\]%
\[
\rho\left(  x,t\right)  =\left(  t^{2}+\sum_{k=1}^{n}x_{k}{}^{2\ell}\right)
^{1/2\ell}\quad.
\]
Thus Theorem \ref{main} applies, and we conclude that the mapping (\ref{map})
is an isomorphism if and only if%
\[
-\frac{n+\ell}{p}<s<n-\frac{n+\ell}{p}\quad.
\]
For the special case of the heat equation,%
\[
\mathcal{L}u=\Delta u-\partial u/\partial t\quad,
\]
this result was proved in \cite{HM}.
\end{example}

\begin{example}
Consider in $\mathbb{R}^{n}$ the operator%
\[
\mathcal{L}u=\sum_{\left\vert \alpha\right\vert =\ell}A_{\alpha}\partial
_{x}^{\alpha}u\quad,
\]
where again each $A_{\alpha}$ is a complex $m\times m$ matrix. The symbol is%
\[
L\left(  x\right)  =\sum_{\left\vert \alpha\right\vert =\ell}A_{\alpha}\left(
ix\right)  ^{\alpha}\quad,
\]
and the semiellipticity requirement that $L\left(  x\right)  $ be invertible
if $x\neq0$ reduces to the usual requirement for ellipticity of the operator.
We have%
\[
\underline{\ell}=\left(  \ell,\ell,\ldots,\ell\right)  \qquad,\qquad
\gamma=\left(  1,1,\ldots,1\right)  \qquad,\qquad\left\Vert \gamma\right\Vert
=n\quad,
\]%
\[
\rho\left(  x\right)  =\left(  \sum_{k=1}^{n}x_{k}{}^{2\ell}\right)
^{1/2\ell}\quad.
\]
Obviously $\rho\left(  x\right)  $ is equivalent to the simpler weight
function $\left\vert x\right\vert $. Theorem \ref{main} applies only if
$\left\Vert \gamma\right\Vert =n>\ell$, in which case the mapping (\ref{map})
is an isomorphism if and only if%
\[
-\frac{n}{p}<s<n-\ell-\frac{n}{p}\quad.
\]

\end{example}

\begin{example}
Let $k$ and $r$ be positive integers, and let $\mathcal{L}$ be the operator%
\begin{equation}
\mathcal{L}u=\sum_{j=0}^{k}\ \sum_{\left\vert \beta\right\vert =jr}%
A_{\beta,k-j}\partial_{x}^{\beta}\partial_{t}^{k-j}u\quad,\label{pparab}%
\end{equation}
where again $u=u\left(  x_{1},\ldots,x_{n},t\right)  $, and each $A_{\beta}$
is a complex constant $m\times m$ matrix. The semiellipticity condition is
that the matrix%
\begin{equation}
\sum_{j=0}^{k}\ \sum_{\left\vert \beta\right\vert =jr}A_{\beta,k-j}\left(
ix\right)  ^{\beta}\left(  it\right)  ^{k-j}\label{ppar}%
\end{equation}
be invertible when $\left(  x,t\right)  \neq0$. The order of this operator is
$\ell=kr$, while%
\[
\underline{\ell}=\left(  kr,\ldots,kr,k\right)  \qquad,\qquad\gamma=\left(
1,\ldots,1,r\right)  \qquad,\qquad\left\Vert \gamma\right\Vert =n+r\quad,
\]%
\[
\rho\left(  x,t\right)  =\left(  t^{2k}+\sum_{j=1}^{n}x_{k}{}^{2kr}\right)
^{1/2kr}\quad.
\]
Our condition $\left\Vert \gamma\right\Vert >\ell$ reduces to $n>\left(
k-1\right)  r$, and the isomorphism condition (\ref{iso}) becomes%
\[
-\frac{n+r}{p}<s<n-\left(  k-1\right)  r-\frac{n+r}{p}\quad.
\]

In the scalar case $m=1$, and with $A_{0,k}=i^{-k}$, (\ref{ppar}) becomes%
\[
f\left(  x,t\right)  :=t^{k}+\sum_{j=1}^{k}\ \sum_{\left\vert \beta\right\vert
=jr}A_{\beta,k-j}\left(  ix\right)  ^{\beta}\left(  it\right)  ^{k-j}\quad.
\]
The operator $\mathcal{L}$ in this case is said to be \textquotedblleft
r-parabolic" \textit{\cite{IT,MIZ} under the additional assumption that
}$\operatorname{Im}\ t\geq\delta>0$ for each root $t$ of $f\left(  x,t\right)
=0$, as $x $ ranges over $\mathbb{R}^{n}$ with $\left\vert x\right\vert =1$.
\end{example}

\end{document}